\documentclass{amsart}
\usepackage{amssymb}
\usepackage{amscd}
\usepackage{amsmath}
\usepackage{amsthm}

\newcommand{\Q}{\mathbb Q}
\newcommand{\Z}{\mathbb Z}
\newcommand{\R}{\mathbb R}
\newcommand{\C}{\mathbb C}

\newcommand{\F}{\mathbb F}

\newcommand{\A}{\mathbb A}
\newcommand{\N}{\mathbb N}

\newcommand{\mQ}{\mathbf{Q}}

\newcommand{\fL}{\mathfrak{L}}
\newcommand{\fO}{\mathfrak{O}}

\newcommand{\tensor}{\otimes}

\newcommand{\fI}{\mathfrak{i}}
\newcommand{\Set}[1]{\left\{\,#1\,\right\}}

\makeatletter
\def\l@section{\@tocline{1}{4pt}{1pc}{}{}}
\def\l@subsection{\@tocline{2}{0pt}{2pc}{5pc}{}}
\makeatother

\begin{document}
\title{On the exceptional zeros of Rankin-Selberg $L$-functions}
\author{Dinakar Ramakrishnan\footnote{Partially supported by the NSF grant DMS-9801328} and Song Wang}
\address
{253-37 Caltech, Pasadena, CA 91125, USA}
\email{dinakar@its.caltech.edu \quad songw@its.caltech.edu}
\maketitle
\pagestyle{myheadings}
\markboth{D.Ramakrishnan and S.Wang}{Exceptional zeros of Rankin-Selberg $L$-functions}

\section*
{\bf Introduction}

\bigskip

In this paper we study the possibility of real zeros near $s=1$ for the Rankin-Selberg $L$-functions 
$L(s, f \times g)$ and $L(s, {\rm sym}^2(g) \times {\rm sym}^2(g))$,
where $f,g$ are newforms, holomorphic or otherwise, on the upper half plane $\mathcal H$, and sym$^2(g)$ denotes the automorphic form on GL$(3)/\Q$ associated to $g$ by Gelbart and Jacquet ([GJ79]). We prove that the set of such zeros of these $L$-functions is the union of the corresponding sets for $L(s, \chi)$ with $\chi$ a quadratic Dirichlet character, which divide them. Such a divisibility does not occur in general, for example when $f, g$ are of level $1$. When $f$ is a Maass form for SL$(2, \Z)$ of Laplacian eigenvalue $\lambda$, this leads to a sharp lower bound, in terms of $\lambda$, for the norm of sym$^2(f)$ on GL$(3)/\Q$, analogous to the well known and oft-used result for the Petersson norm of $f$ proved in [HL94] 
and [GHLL94]. As a consequence of our result on $L(s, {\rm sym}^2(g) \times {\rm sym}^2(g))$ one gets a good upper 
bound for the {\it spectrally normalized} first coefficient $a(1,1)$ of sym$^2(g)$. (In the artihmetic normalization, $a(1,1)$ would be $1$.) In a different direction, we are able to show that the symmetric sixth and eighth power $L$-functions of modular forms $f$ with trivial character (Haupttypus) are holomorphic in $(1 - \frac{c}{\log M}, 1)$, where $M$ is the {\it thickened conductor} (see section 1) and $c$ a universal, positive, effective constant; by a recent theorem of Kim and Shahidi
([KSh2001]), one knows that these $L$-functions are invertible in $\Re(s) \geq 1$ except possibly for a pole at $s=1$. If $f$ runs over holomorphic newforms of a {\it fixed} weight (resp. level), for example, the thickened conductor $M$ is essentially the level (resp. weight). We will in general work over arbitrary number fields and use the adelic language.

\medskip

First some preliminaries. Suppose $D(s)$ is any Dirichlet series given as an Euler product in $\{\Re(s) > 1\}$, which admits a meromorphic continuation to the whole $s$-plane with no pole outside $s=1$, together with a functional equation relating $s$ to $1-s$ after adding suitable archimedean factors. By an {\it exceptional zero}, or a {\it Siegel zero}, or perhaps more appropriately (cf. [IwS2000]) a {\it Landau-Siegel zero}, of $D(s)$, one means a real zero $s = \beta$ of $D(s)$ which is close to $s=1$. More precisely, such a zero will lie in $(1-\frac{C}{\log M}, 1)$, where $C$ is an effective, universal constant $> 0$ (see section 1). The {\it Grand Riemann Hypothesis} (GRH) would imply that there should be no such an exceptional zero, but it is of course quite hard to verify.

\medskip

It was shown in [HRa95] that for any number field $F$, the $L$-function $L(s, \pi)$ of a cusp form $\pi$ in GL$(2)/F$ admits no Landau-Siegel zero. In the special case when $\pi$ is {\it dihedral}, i.e., associated to a character $\chi$ of a quadratic extension $K$ of $F$, $L(s,\pi)$ is simply the abelian $L$-function $L(s,\chi)$ considered by Hecke, and if $\theta$ is the non-trivial automorphism of $K/F$, the cuspidality of $\pi$ forces $\chi$ to be distinct from $\chi \circ \theta$. We will say that $\pi$ is of type $(K/F,\chi)$ in this case. 

It was also shown in [HRa95] that for any $n > 1$, the standard $L$-series $L(s, \pi)$ of cusp forms $\pi$ on GL$(n)/F$ admit no Landau-Siegel zero {\it if} one assumes Langlands's {\it principle of functoriality}, in particular the existence of the {\it automorphic tensor product}. An analogous, but slightly more complicated, statement can be made for general Rankin-Selberg $L$-series on GL$(n) \times $GL$(m)$, but assuming the full force of functoriality is but a distant dream at the moment, though it is highly instructive to be aware of what it entails. So it becomes an interesting problem to know how much one can {\it unconditionally} prove by making use of available instances of functoriality; the method has to deviate some from that given in [HRa95]. This is what we carry out here for $n = m \leq 3$. 

Roughly speaking, the main point is to find a suitable positive Dirichlet series $D(s)$ which is divisible by the $L(s)$ of interest to a degree $k$, say, which is (strictly) {\it larger} than the order of pole of $D(s)$ at $s=1$. If there is anything creative here, at all, it is in the proper choice of $D(s)$ and then in the verification of the holomorphy of $D(s)/L(s)^k$, at least in a real interval $(t,1)$ for a fixed $t < 1$. It should be noted, however, that this approach fails to give anything significant for $L$-functions of quadratic characters; for two very interesting, and completely different, approaches for this crucial case see [IwS2000] and [GS2000].

\medskip

Now fix a number field $F$ and consider the Rankin-Selberg $L$-function $L(s, \pi \times \pi')$ associated to a pair $(\pi, \pi')$ of cusp forms on GL$(2)/F$. Denote by $\omega$, resp. $\omega'$, the central character of $\pi$, resp. $\pi'$. Our first main result is the following

\medskip

\noindent{\bf Theorem A} \, \it Let $\pi, \pi'$ be cuspidal automorphic representations of GL$(2, \A_F)$. Then $L(s, \pi \times \pi')$ admits no Landau-Siegel zero except possibly in the following cases:
\begin{enumerate}
\item[{(i)}] $\pi$ is non-dihedral and $\pi' \simeq \pi \otimes \mu$ with $\omega\mu$ of order $\leq 2$;
\item[{(ii)}] $\pi$, resp. $\pi'$, is dihedral of type $(K,\chi)$, resp. $(K',\chi')$, with $K' = K$ and $\chi'\chi$ or $\chi'(\chi \circ \theta)$ of order $\leq 2$.
\end{enumerate}
In case (i), resp. (ii), the exceptional zeros of $L(s, \pi \times \pi')$ are the same as those of $L(s,\omega\mu)$, resp. $L(s,\chi'\chi)L(s,\chi'(\chi \circ \theta))$. 
In case (ii), if $\chi'\chi$ or $\chi'(\chi \circ \theta)$ is trivial, then the exceptional zeros are the same as those of $\zeta_K(s)$. 
In either case, there is at most one exceptional zero.
\rm

\medskip

For the vast majority of cases not satisfying (i) or (ii), $L(s, \pi \times \pi')$ has no exceptional zero. In particular, if $\pi_0, \pi'_0$ are fixed, non-dihedral cusp forms on GL$(2)/F$ which are not twist equivalent to each other, there exists an effective constant $c > 0$ such that the family $L(s, \pi_0 \times (\pi'_0 \otimes \chi))$, with $\chi$ running over quadratic characters of conductor $q$ prime to the levels of $\pi_0, \pi'_0$, admits {\it no} real zero $\beta$ with $\beta \in \left(1 - \frac{c}{{\rm log}q}\right)$. 

\medskip

In case (i) we have
$$
L(s, \pi \times \pi') \, = \, L(s,\omega\mu)L(s, sym^2(\pi) \otimes \mu).
$$
The non-existence of Landau-Siegel zeros for $L(s, sym^2(\pi))$ (for $\pi$ non-dihedral) has been known for a while by the important work of Goldfeld, Hoffstein, Lieman and Lockart ([GHLL94]). For general $\mu$, the non-existence for $L(s, sym^2(\pi) \otimes \mu)$ is known by [Ba97], following an earlier reduction step given in section 6 of [HRa95]. So our Theorem is not new in this case. 

\medskip

For any cusp form $\pi$ on $GL(2)/F$, let $L(s, \pi; {\rm sym}^n)$ denote, for every $n \geq 1$, the symmetric $n$-th power $L$-function of $\pi$ (see section 1 for a definition). It is expected that there is an automorphic form sym$^n(\pi)$ on GL$(n+1)/F$ whose standard $L$-function coincides with $L(s, \pi; {\rm sym}^n)$. This is classical ([GJ79]) for $n=2$ and a major breakthrough has been made recently for $n=3$ ([KSh2000]) and for $n=4$ ([K2000]). 
The proof of Theorem A uses the result for $n=3$ as well as the construction of the first author ([Ra2000]) of the Rankin-Selberg product of  pairs $(\pi, \pi')$ of forms on GL$(2)/F$ as an automorphic form $\pi \boxtimes \pi'$ on GL$(4)/F$.

Recall that $\pi$ is {\it dihedral} iff it admits a self-twist by a non-trivial, necessarily quadratic, character. One says that it is {\it tetrahedral}, resp. {\it octahedral}, iff sym$^2(\pi)$, resp. sym$^3(\pi)$, is cuspidal and admits a non-trivial self-twist by a cubic, resp. quadratic character. It is well known that sym$^2(\pi)$ is cuspidal ([GJ79]) iff $\pi$ is not dihedral. It has been shown in [KSh2000], resp. [KSh2001], that sym$^3(\pi)$, resp. sym$^4(\pi)$, is cuspidal iff $\pi$ is not dihedral or tetrahedral, resp. not dihedral, tetrahedral or octahedral.
We will henceforth say that a cusp form $\pi$ on $GL(2)/F$ is of {\it solvable polydedral type} 
if it is either dihedral or tetrahedral or octahedral. Our second main result is the following:

\medskip

\noindent{\bf Theorem B} \, \it Let $\pi$ be a self-dual cuspidal automorphic representation of GL$(2, \A_F)$. Then the set of Landau-Siegel zeros of $L(s, {\rm sym}^2(\pi) \times {\rm sym}^2(\pi))$ is the union of the sets of Landau-Siegel zeros of abelian $L$-functions of the form $L(s, \chi)$, $\chi^2 = 1$, which divide $L(s, {\rm sym}^2(\pi) \times {\rm sym}^2(\pi))$, if any. Moreover, if $\pi$ is not of solvable polyhedral type and if $\zeta_F(s)$ has no exceptional zero (for example when $F = \Q$), then there is no exceptional zero.
\rm

\medskip

This theorem holds also for any cuspidal $\pi$ on GL$(2)/F$ which is a unitary character twist of a selfdual representation. 

\medskip

The following corollary of Theorem B gives precise bounds for the Petersson norm on $GL(3)/\mQ$
of the symmetric square of a Maass wave form $g$ of level $1$.
In Section 6 (see Definition 6.3), we will define a suitably normalized function ${\rm sym}^{2} (g)$ 
spanning the space of
the Gelbart-Jacquet square lift to $GL(3)/\mQ$ of the cuspidal automorphic
representation $\pi$ generated by $g$.

\medskip

\noindent{\bf Corollary C} \,\it
Let $g$ be a Maass form on the
upper half plane $\mathcal H$, relative to SL$(2, \Z)$, 
of weight zero and Laplacian eigenvalue $\lambda$, which
is also an eigenfunction of Hecke operators.
Then for each $\varepsilon > 0$,
	\[
		\frac{1}{\log (\lambda + 1)}
		 << \langle\,{\rm sym}^{2} (g), {\rm sym}^{2} (g)\,\rangle
		<<_{\varepsilon} (\lambda + 1)^{\varepsilon},
	\]
	where ${\rm sym}^2(g)$ is spectrally mormalized as in Theorem B.

Moreover, if $\{a(m,n)\}$ denotes the collection of Fourier coefficients of the spectrally normalized function
${\rm sym}^2(g)/\vert\vert{\rm sym}^2(g)\vert\vert$, we have
$$
\vert a(1,1)\vert \, << \, \log (\lambda + 1).
$$
\rm

\medskip

Our proof of Theorem B will establish on the way that the symmetric $4$-th power $L$-function of any self-dual cusp form $\pi$ on GL$(2)$, not of solvable polyhedral type, admits no Landau-Siegel zero. A slew of interesting results for the $L(s, \pi; {\rm sym}^n)$ for {\it $n$ up to $9$} have been established recently by H.~Kim and F.~Shahidi in [KSh2001], proving in particular the meromorphic continuation, functional equation and holomorphy in $\Re(s) \geq 1$, except for a possible pole at $s=1$. It may be of some interest to know, for $n > 4$, how far to the left of $s=1$ can the holomorphy assertion can be extended. A consequence of our work is the following tiny, but apparently non-trivial, extension to the left of $s=1$ for $n=6, 8$.

\medskip

\noindent{\bf Theorem D} \, \it Let $F$ be a number field and $\pi$ a self-dual cusp form on GL$(2)/F$ of thickened conductor $M$. Then there exists a universal, effective constant $c > 0$ such that $L(s, \pi; {\rm sym}^6)$ and $L(s, \pi; {\rm sym}^8)$ have no pole in the real interval $(1 - \frac{c}{\log M}, 1)$.
\rm

\medskip

Now we will say a few words about the proofs.

Regarding Theorem A, suppose we are in the main case, i.e., neither $\pi$ nor $\pi'$ is dihedral and also $\pi'$ is not isomorphic to $\pi \otimes \mu$ for any character $\mu$.
Under these hypotheses, $\pi \boxtimes \pi'$ is cuspidal on GL$(4)/F$ by [Ra2000]. When it is not self-dual, there is a simple argument (see section 3 of [HRa95]) to deduce the non-existence of a Siegel zero. So we may assume that $\pi \boxtimes \pi'$ is self-dual, which implies that the central characters $\omega, \omega'$ of $\pi, \pi'$ are inverses of each other. For simplicity assume for the moment that $\omega, \omega'$ are trivial. (For a full treatment of the general case, see section 4.) 
Then the key point is to appeal to the following identity of $L$-functions
\begin{align}
L(s, \Pi \times \Pi) = &\zeta_F(s) L(s, {\rm sym}^2(\pi))^2 L(s, \pi \times \pi)^4 L(s, {\rm sym}^3(\pi) \times \pi')^2 \notag \\
&L(s, {\rm sym}^2(\pi) \times {\rm sym}^2(\pi)) L(s, (\pi \boxtimes \pi') \times (\pi \boxtimes \pi')).
\notag
\end{align}
where $\Pi$ is an isobaric automorphic form on GL$(8)/F$ defined by 
$$
\Pi \, = \,  1 \boxplus \pi \boxtimes \pi' \boxplus
{\rm sym}^{2} ({\pi}),
$$
where $1$ denotes the trivial automorphic form on GL$(1)/F$ and $\boxplus$ the Langlands sum operation on automorphic forms defined by his theory of Eisenstein series ([La79]), proved to be well defined by the work of Jacquet and Shalika ([JS81]). The degree $64$ Rankin-Selberg $L$-function $L(s, \Pi \times \Pi)$ has the standard analytic properties and defines, in $\Re(s) > 1$, a Dirichlet series with non-negative coefficieents. Moreover, it has a pole of order $3$ at $s=1$, and since since $L(s, \pi \times \pi')$ occurs in its factorization to a power larger than $3$, a standard lemma (see Lemma 1.7) precludes the latter from having any Landau-Siegel zero. We also need to show that the ratio $L(s, \Pi \times \Pi)/L(s, \pi \times \pi')^4$ is holomorphic, for which we appeal to the automorphy of sym$^3(\pi)$ ([KSh2000]).

\medskip

The proof of Theorem B involves a further wrinkle, and uses in addition the automorphy of sym$^4(\pi)$ ([K2000]), as well as the works of Bump-Friedberg ([BuG92]) on the symmetric square $L$-functions of GL$(n)$. The well known identity
$$
L(s, {\rm sym}^2(\pi) \times {\rm sym}^2) \, \zeta_F(s) L(s, {\rm sym}^2(\pi)) L(s, {\rm sym}^4(\pi))
$$ 
reduces the problem to studying the Landau-Siegel zeros of $L(s, {\rm sym}^4(\pi))$. We show in section 5 (see Theorem B$^\prime$) that for $\pi$ {\it not} of solvable polyhedral type, $L(s, \pi; {\rm sym}^4)$ has no exceptional zero. To do this we set
$$
\Pi : \, = \, 1 \boxtimes {\rm sym}^2(\pi) \boxplus {\rm sym}^4(\pi),
$$
and consider
\begin{align}
L(s, \Pi \times \Pi) = &\zeta_F(s) L(s, {\rm sym}^2(\pi))^2 L(s, {\rm sym}^4(\pi))^4 L(s, \pi; {\rm sym}^6)^2
\notag \\
&L(s, {\rm sym}^2(\pi) \times {\rm sym}^2(\pi)) L(s, {\rm sym}^4(\pi) \times {\rm sym}^4(\pi)).
\notag
\end{align}
Since $L(s, \Pi_f \times \Pi_f)$ defines a  Dirichlet series in $\Re(s) > 1$ with non-negative coefficients and has a pole of order  $3$ at $s=1$, things appear to be in good shape, {\it till} we realize that one does not yet know how to prove that $L(s, \pi; {\rm sym}^6)$ is holomorphic in any real interval $(1-t, 1)$ for a fixed $t < 1$. But luckily we also have the factorization
$$
L(s, {\rm sym}^3(\pi); {\rm sym}^2) \, = \, L(s, \pi; {\rm sym}^6) L(s, {\rm sym}^2(\pi)),
$$
which allows us, after exercising some care about the bad factors, to prove the holomorphy of the ratio $L(s, \Pi \times \Pi)/L(s, {\rm sym}^4(\pi))^4$ in $(1/2, 1)$. 
For details we refer to section 5.

\medskip

In order to prove Corollary C, we begin by deducing a precise relationship between
the Peterson norm of a suitably mormalized ${\rm sym}^2(g)$ (see section 6) and 
$L(s, {\rm sym}^{2} (g) \times {\rm sym}^{2} (g))$ by
using the results on its integral representation due to Jacquet, Piatetski-Shapiro and Shalika ([JPSS83], [JS90]), 
and an exact formula of E. Stade for pairs of spherical representations of GL$(3, \R)$ ([St93], [St2001]). More precisely, $(\,{\rm sym}^{2}(g), {\rm sym}^{2}(g)\,)$
differs $\text{\rm Res}_{s = 1} L(s, {\rm sym}^{2}(g) \times {\rm sym}^{2}(g))$ by 
a constant factor coming from the residue of an Eisenstein series. (There should also be
a similar result when $f$ is a holomorphic newform for SL$(2, \Z)$, but 
at this point one does
not appear to know enough about the archimedean zeta integral on GL$(3) \times $GL$(3)$ to achieve this.)
Finally, the nonexistence of Landau-Siegel zeros for
$L(s, {\rm sym}^{2}(g) \times {\rm sym}^{2}(g))$ allows us to bound from below its residue at $s=1$, and this in turn gives us, when sym$^2(g)$ is replaced by its spectral normalization, the desired bound on the {\it first Fourier coefficient} $a(1,1)$.
For details, see Section 6.

\medskip

To prove Theorem D we appeal to the factorization of $L(s, {\rm sym}^4(\pi); \Lambda^2)$ above as well as to the identity 
$$
L(s, {\rm sym}^4(\pi); {\rm sym}^2) \, = \, L(s, \pi; {\rm sym}^8) L(s, {\rm sym}^4(\pi)) \zeta_F(s).
$$
For any cuspidal automorphic representation $\Pi$ of GL$(n, \A_F)$, if $S$ denotes the union of the archimedean places of $F$ with the set of finite places where $\Pi$ is ramified, one knows by Bump-Ginzburg ([BuG92]) that the incomplete $L$-function $L^S(s, \Pi; {\rm sym}^2)$, defined in a right half plane by the Euler product over places outside $S$, is holomorphic in $(1/2, 1)$. Applying this with $n=4, \Pi = {\rm sym}^4(\pi)$, and carefully taking care of the factors at $S$, we deuce the holomorphy of the syemmetric square $L$-function of sym$^4(\pi)$. Now the knowledge gained from the proof of Theorem B on Landau-Siegel zeros of $L(s, {\rm sym}^4(\pi))$ allows us  (see section 7) to prove Theorem D.

\medskip

We would like to thank D.~Bump, W.~Duke, H.~Jacquet, H.~Kim, W.~Luo, S.~Miller, F.~Shahidi and E.~Stade for useful conversations and/or correspondence. Clearly this paper depnds on the ideas and results of the articles [HL94], [GHLL94], [HRa95], [Ra2000], [KSh2000,1] and [K2000]. The first author would also like to thank the NSF for support through the grants DMS-9801328 and DMS-0100372. The subject matter of this paper formed a portion of the first author's {\it Schur Lecture} at the University of Tel Aviv in March 2001, and he would like to thank J.~Bernstein and S.~Gelbart for their interest.

\vskip 0.2in

\section{\bf Preliminaries on Landau-Siegel zeros}

\bigskip

For every $m \geq 1$, 
let ${\mathcal D}_m$ denote the class of {\it
Dirichlet series} $L(s)= \sum\limits_{n \geq 1} \frac{a_n}
{n^s}$, absolutely convergent in $Re(s) > 1$ with an Euler
product $\prod_p P_p(p^{-s})^{-1}$ of degree $m$
there, extending to whole $s-plane$ as a meromorphic function of
bounded order, in fact with no poles anywhere except at $s=1$,
and satisfying (relative to another Dirichlet series
$L^{\vee}(s)$ in ${\mathcal D}_m$) a functional equation of the form   
$$
L_\infty (s) \, L(s) = WN^{1/2 - s} L_\infty ^{\vee} (1 -
s) L^{\vee} (1-s),
\leqno(1.1)
$$
where $W \in \C^\ast$, $N \in \N$, 
and the {\it archimedean factor} $L_\infty (s)$ is 
$$
\pi^{-sm/2} \prod\limits_{j=1}^m \Gamma\left(\frac{s + b_j}{2}\right),
$$
where $(b_j) \, \in \, \C^m$.

\medskip

Put ${\mathcal D} = \bigcup\limits_{m \geq 1} {\mathcal D}_m$.
It might be useful to compare this with the definition of the Selberg class ([Mu94]). 
In the latter one requires in addition a Ramanujan type bound on the coefficients, but allows more complicated Gamma factors.

\medskip

One says that $L(s)$ is {\it self-dual} if $L(s)=L^{\vee}(s)$, in which
case $W = \{\pm 1 \}$.

\medskip

Important examples are $\zeta(s)$, Dirichlet and Hecke $L$-functions, the $L$-functions
of {\it holomorphic newforms} $g$ of weight $k \geq 1$ and level $N$, normalized to be
$$
L(s) = L (s + \frac{k-1}{2}, g), \,\, L_\infty(s) = \pi^{-s}
\Gamma\left(\frac{s+(k-1)/2}{2}\right) \Gamma\left(\frac{s+(k+1)/2}{2}\right),
\leqno(1.2)
$$
$L$-functions of {\it cuspidal Maass forms} $\phi$ of level $N$ which are 
eigenfunctions of Hecke
operators:
$$
L(s) = L(s, \phi), \,\, L_\infty(s)=\pi^{-s} \Gamma\left(\frac{s + \delta +w}{2}\right) 
\Gamma\left(\frac{s+\delta - w}{2}\right), \,\, \delta \in \{0,1 \},
\leqno(1.3)
$$
and the {\it Rankin-Selberg $L$-functions} $L(s, f \times g)$, where $f, g$ are cusp forms of 
holomorphic or Maass type. See the next section for their definition and generalizations. 

\medskip

Call an $L(s) \in \mathcal D$ {\it primitive} if it cannot be factored as $L_1(s)L_2(s)$ with $L_1(s), L_2(s)$ non-scalar in $\mathcal D$. The Dirichlet and Hecke $L$-functions, as well as those of cusp forms on the upper half plane are primitive.

\medskip

\noindent{\bf Conjecture I} \, \it Every $L(s) \, \epsilon \, {\mathcal D}_m$
is {\it quasi-automorphic}, ie, there exists an automorphic form
$\pi$ on GL$(m)/\Q$ such that $L_p(s)\, = \, L(s, \pi_p)$ for almost all $p$.
Moreover, $L(s)$ is primitive iff $\pi$ is cuspidal.
\rm

\medskip

This is compatible with Langlands philosophy ([La70]) and with the 
conjecture of Cogdell and Piatetski-Shapiro ([CoPS94]); Cogdell has remarked to us 
recently that Piatetski-Shapiro has also formulated (unpublished) a similar conjecture 
involving analogs of $\mathcal D$.
Note such a thing cannot be formulated over number fields of degree $> 1$ as one can
permute the Euler factors lying over any given rational prime.  
Also, there exists an example of
Patterson over function fields $F$ over a finite field $\F_q$ satisfying
analogous conditions, but with zeros on the lines $Re(s) = 1/4$ and $Re(s) = 3/4$. 
One problem in characteristic $p$ is that there is no
minimal global field $F$ such as $\Q$. It is still an interesting open problem to know 
if one can define a good notion of "primitivity" over function fields.

\medskip

\noindent{\bf Conjecture II} \, \it For any $L(s) \, \in \,
{\mathcal D}$, if it has a pole of order $r$ at $s=1$, then
$\zeta(s)^r | L(s)$, ie., $L(s)=\zeta(s)^r L_1(s)$, with $L_1(s)
\, \epsilon \, {\mathcal D}$. 
\rm

\medskip

This is compatible with the conjectures of Selberg, Tate and
Langlands. 

\bigskip

\noindent{\bf Definition 1.4} \, \it Let $L(s) \,
\epsilon \, {\mathcal D}_m$ with $L_{\infty}(s) = \pi^{-ms} 
\mathop\Pi_{j=1}^m
\Gamma\left(\frac{s+ b_j}{2}\right)$.  Define its {\it thickened conductor}
to be
$$
M = N (2+ \Lambda)
$$
where
$$
\Lambda \, = \, \mathop\sum_{j=1}^m |b_j|).
$$
\rm

\medskip

\noindent{\bf Definition 1.5} \, \it Let $c>0$.  Then we say that $L(s)$
has a Landau-Siegel zero relative to $c$ if $L(\beta)=0$ for
some $\beta \, \in \, (1- \frac{c}{\log M}, 1)$.
\rm

\medskip

\noindent{\bf Definition 1.6} \it Let ${\mathcal F}$ be a family, by wich we will mean a
class of $L$-functions in ${\mathcal D}$ with $M
\rightarrow \infty$ in ${\mathcal F}$.  We say that ${\mathcal F}$
admits no {\it Landau-Siegel Zero} if there exists an effective
constant $c>0$ such that no $L(s)$ in ${\mathcal F}$ has a zero in
$(1- \frac{c}{\log M}, 1)$. 
\rm

\medskip

The general expectation is framed by the following

\medskip

\noindent{\bf Conjecture III} \, \it Let $\mathcal F$ be a family in $\mathcal
D$.  Then $\mathcal F$ admits no Landau-Siegel zero. 
\rm

\medskip

One reason for interest in this is that the lack of such a zero
implies a good lower bound for $L(s)$ at $s=1$. Note that in view of Conjecture I, the 
{\it Grand Riemann Hypothesis}, shortened as GRH, implies that all the non-trivial zeros of any
$L(s)$ in $\mathcal D$ lie on the critical line, hence it implies Conjecture III. 
Of course the GRH is but a distant goal at the moment, and it is hopefully of interest to 
verify Conjecture III for various families.

The $L$-functions of pure motives over $\Q$, in particular those associated to the cohomology 
of smooth projective varieties $X/\Q$, are expected to be automorphic and hence should belong to 
$\mathcal D$. For these $L$-functions, when they are of {\it even} Frobenius weight, the values at 
$s=1$ have arithmetic significance by the general Bloch-Kato conjectures, and so the question 
of non-existence of Landau-Siegel zeros is helpful to understand from a purely arithmetical point of view. 

\medskip

We will need the following useful (and well known) fact:

\medskip

\noindent{\bf Lemma 1.7} \, \it Let $L(s) \, \epsilon \, {\mathcal D}_m$ be a positive
Dirichlet series having a pole of order $r \geq 1$ at $s=1$,
with $L'(s)/L(s) <0$ for real $s$ in $(1, \infty)$.  Then there
exists an
effective constant $C>0$, depending only on $m$ and $r$, such
that $L(s)$ has most $r$ real zeros in $(1 - \frac{C}{\log M}, 1)$.
\rm

\medskip

For a more relaxed discussion of these matters, see the expository article [Ra99] on the Landau-Siegel 
zeros, as well as the articles [GHLL94] and [HRa95].

\vskip 0.2in

\section{\bf Preliminaries on automorphic $L$-functions}

\bigskip

Fix a number field $F$ with ring of integers $\fO_F$, discriminant $D_F$, and adele ring $\A_F = F_\infty \times \A_{F,f}$, where $F_\infty$ is the product of the archimedean completions of $F$ and the ring $\A_{F,f}$ of finite adeles is a restricted direct product of the completions $F_v$ over non-archimedean places $v$. For each finite $v$, let $\fO_v$ denote the ring of integers of $F_v$. When $F = \Q$, $F_\infty \simeq \R$ and $\A_{F,f} \simeq \hat \Z \otimes \Q$, where $\hat \Z$ is the inverse limit of $\{\Z/m \vert m \geq 1\}$ and is non-canonically isomorphic to $\prod_p \Z_p$. 

Recall that a cuspidal automorphic representation $\pi$ of GL$(n, \A_F)$ is among other things admissible, i.e., a restricted tensor product $\otimes_v \pi_v \simeq \pi_\infty \otimes \pi_f$, where $v$ run over all the places of $F$ and $\pi_v$ is, for almost all finite $v$, unramified, i.e., its space admits a vector invariant under GL$(n, \fO_v)$. Given a partition of $n$ as $\sum_{j =1}^r n_j$ with each $n_j \geq 1$, and cuspidal automorphic representations $\pi_1, \ldots, \pi_r$ of GL$(n_1, \A_F), \ldots $GL$(n_r,\A_F)$, Langlands's theory of Eisenstein series constructs a so called  {\it isobaric} automorphic representation ([La79]) $\pi$ of GL$(n, \A_F)$, which is unique by the work of Jacquet-Shalika ([JS81]), written as
$$
\pi: = \, \boxplus_{i=1}^r \pi_i,
\leqno(2.1)
$$
with the property that its standard degree $n$ $L$-function $L(s, \pi)$ is the product $\prod_{i=1}^r L(s, \pi_i)$. 
Write
$$
L(s, \pi_{\infty}) \, = \, \pi^{-dns/2}\prod_{j=1}^{dn} \Gamma(\frac{s+b_j(\pi)}{2}),
\leqno(2.2)
$$
where $d = [F:\Q]$ and the $b_j(\pi)$ are complex numbers depending only on $\pi_\infty$.

Now consider a pair of isobaric automorphic representations $\pi, \pi'$ of ${\rm GL}(n, \A_F)$, GL$(m, \A_F)$,
respectively. The associated Rankin-Selberg $L$-function is given as an Euler product of degree $nm$:
$$
L(s, \pi \times \pi') \, = \, \prod_v L(s, \pi_v \times \pi'_v),
\leqno(2.3)
$$
convergent in a right half plane, with its {\it finite part}, namely $L(s, \pi_f \times \pi'_f)$, defining a Dirichlet series. When $m=1$ and $\pi'$ is the trivial representation $1$, this $L$-function agrees with the standard $L$-function. 
There are two distinct methods for defining these $L$-functions, the first using the {\it gcd}s of
integral representations, due to Jacquet, Piaietski-Shapiro and Shalika ([JPSS83]), and the second via the constant terms of Eisenstein series on larger groups, due to Langlands and Shahidi ([Sh88, 90]); see also [MW89]. 
The fact that they give the same
$L$-functions is non-trivial but true. These $L$-functions 
also admit a meromorphic continuation to the whole $s$-plane with no poles except possibly at 
$1-s_0$ and $s_0$ for a unique $s_0 \in i\R$; such a pole occurs iff $\pi$ and 
$\pi' \otimes \vert.\vert^{s_0}$ are contragredients of each other. One also has the functional equation
$$
L(s, \pi \times \pi') \, = \, 
\varepsilon(s, \pi \times \pi')W(\pi \times \pi') L(1-s, \overline \pi \times \overline \pi')
\leqno(2.4)
$$
where 
$$
\varepsilon(s, \pi \times \pi') \, = \, (d_F^{nm}N(\pi \times \pi'))^{\frac{1}{2}-s},
$$
which is an invertible holomorphic function. Here $N(\pi \times \pi')$ is the {\it conductor},  
and $W(\pi \times \pi') \in \C^\ast$ the {\it root number}, of the pair $(\pi, \pi')$.

The following was proved in [HRa95] (Lemma $a$ of section 2):

\medskip

\noindent{\bf Lemma 2.5} \, \it For any unitary, isobaric automorphic representation $\pi$ of GL$(n, \A_F)$, the
Dirichlet series defined by $L(s, \pi_f \times \overline \pi_f)$ has non-negative coefficients.
Moreover, the logarithmic derivative $L'(s, \pi_f \times \overline \pi_f)/L(s, \pi_f \times \overline \pi_f)$
is negative for real $s$ in $(1, \infty)$.
\rm

\bigskip

The local Langlands correspondence for GL$(n)$, proved by Harris-Taylor ([HaT2000]) and Henniart ([He2000]) in the non-archimedean case (and proved long ago by Langlands in the archimedean case), gives a bijection at any place $v$, preserving the $L$- and $\varepsilon$-factors of pairs, between irreducible admissible representations $\pi_v$ of GL$(n, F_v)$ and $n$-dimensional representations $\sigma_v = \sigma(\pi_v)$ of the extended (resp. usual) Weil group $W_{F_v}': = W_{F_v} \times $SL$(2, \C)$ (resp. $W_{F_v}$) in the $p$-adic (resp. archimedean) case. This gives in particular the identity at any finite $v$:
$$
N(\pi_v \times \pi'_v) \, = \, N(\sigma(\pi_v) \otimes \sigma(\pi'_v)),
\leqno(2.7)
$$
where for any representation $\tau$ of $W_{F_v}'$, $N(\tau)$ denotes the usual Artin conductor. A consequence of this  is the sharp bound:
$$
	{M(\pi)}^{- n'} {M(\pi')}^{- n}
	\le M(\pi \times \pi')
	\le {M(\pi)}^{n'} {M(\pi')}^{n}.
\leqno(2.8)
$$
In fact we do not need the full force of this, and the weaker bound proved in Lemma $b$, section 2 of[HRa95], where the exponents were polynomially dependent on $n, n'$, is actually sufficient for our purposes.

\medskip

Combining all this information with Lemma 1.7 we get 

\medskip

\noindent{\bf Proposition 2.9} \, \it Let $\pi$ be
	an isobaric automorphic representation of $GL(n, \A_{F})$
	with $L(s, \pi \times \bar{\pi})$ having a pole of order
	$r \ge 1$ at $s = 1$. Then there is an effective constant $c \ge 0$
	depending on $n$ and $r$, such that $L(s, \pi \times \overline \pi)$
	has at most $r$ real zeros in the interval
$$
		J : \, = \, \{s \in \C \vert 1 - c / \log M (\pi \times \overline \pi) < \Re(s) < 1\}.
$$
	Furthermore, if $L(s, \pi \times \pi') = {L_{1} (s)}^{k} L_{2} (s)$
	for some nice $L$--series $L_{1} (s)$ and $L_{2} (s)$
	with $k > r$ and $L_{2} (s)$ holomorphic in $(t,1)$ for some fixed $t \in (0,1)$,
	then $L_{1} (s)$ has no zeros in $J$.
\rm

\medskip

This provides a very useful criterion
to prove the nonexistence of Landau-Siegel zeros in some cases. By the definition
of the conductor of the $L$--series, if we know that the logarithm
of the conductor of $L_{2} (s)$ does not exceed some multiple of
the logarithm of the conductor of $L_{1} (s)$, with the 
constant depending only on the degrees of those $L$--series and $k$, then we can conclude
that the logarithm of the conductor of $L_{1} (s)$ is bounded above and below
by some multiples of the logarithm of
the conductor of $L(s, \pi \times \pi')$, which will then imply that $L_{1} (s)$
has no Landau-Siegel zero.

\bigskip

Given any isobaric automorphic representation $\pi$ of GL$(n, \A_F)$, a finite dimensional
$\C$-representation $r$ of (the connected dual group) GL$(n, \C)$, and a character $\mu$ of $W_F$,  
we can define the associated automorphic $L$-function by
$$
L(s, \pi; r \otimes \mu) \, = \, \prod\limits_v \, L(s, r(\sigma(\pi_v)) \otimes \mu_v),
\leqno(2.11)
$$
and
$$
\varepsilon(s, \pi; r \otimes \mu) \, = \, \prod\limits_v \, \varepsilon(s, r(\sigma(\pi_v)) \otimes \mu_v),
$$
where $v$ runs over all the places of $F$, $\pi_v \to \sigma(\pi_v)$ the arrow giving the local Langlands correspondence for GL$(n)/F_v$, and local factors are those attached to representations of the (extended) Weil group ([De73]). (To be precise, in the treatment of the non-archimedean case in [De73], Deligne uses the Weil-Deligne group $WD_{F_v}$, but it is not difficult to see how its {\it representations} are in bijection with those of $W'_{F_v}$. Also, the local $\varepsilon$-factors depend on the choice of a non-trivial additive character and the Haar measure, but we suppress this in our notation.) Originally, Langlands gave a purely automorphic definition of the local factors at almost all places, but now, thanks to [HaT2000] and [He2000], we can do better.

We can also define {\it higher analogs of the Rankin-Selberg $L$-functions} and set, for any pair $(\pi, \pi')$ of isobaric automorphic forms on $({\rm GL}(n), {\rm GL}(m))/F$, a pair $(r, r')$ of finite dimensional $\C$-representations of ${\rm GL}(n, \C), {\rm GL}(m, \C))$, and a character $\mu$ of $W_F$,
$$
L(s, \pi \times \pi'; r \otimes r' \otimes \mu) \, = \, \prod\limits_v \, L(s, r(\sigma(\pi_v)) \otimes r'(\sigma(\pi'_v)) \otimes \mu_v),
\leqno(2.11)
$$
and
$$
\varepsilon(s, \pi \times \pi'; r \otimes r' \otimes \mu) \, = \, \prod\limits_v \, \varepsilon(s, r(\sigma(\pi_v)) \otimes r'(\sigma(\pi'_v)) \otimes \mu_v).
$$

When $m=1$, $\pi' \simeq 1$ and $r' \simeq 1$, $L(s, \pi \times \pi'; r \otimes r' \otimes \mu)$ coincides with $L(s, \pi; r \otimes \mu)$.

\medskip

For each $j \geq 1$, let sym$^j$ denote the symmetric $j$-th power of the standard representation of GL$(n, \C)$. The definition (2.11) above  gives in particular the families of automorphic $L$-functions $L(s, \pi; {\rm sym}^j \otimes \mu)$ and $L(s, \pi \times \pi'; {\rm sym}^j \otimes {\rm sym}^k \otimes \mu)$ for isobaric automorphic representations $\pi, \pi'$ of GL$(2, \A_F)$ and idele class character $\mu$, which we may, and we will, identify (via class field theory) with a character, again denoted by $\mu$, of $W_F$. One calls $L(s, \pi; {\rm sym}^j)$ {\it the symmetric $j$-th power $L$-function}. 

\vskip 0.2in

\section{Some useful instances of functoriality}

\bigskip

Here we summarize certain known instances, which we will need, of functorial transfer of automorphic forms from one group to another.

\medskip

Let $\pi, \pi'$ be cuspidal automorphic representations of GL$(n, \A_F)$, GL$(m, \A_F)$, and let $r, r'$ be $\C$-representations of GL$(n, \C)$, GL$(m, \C)$ of dimension $d, d'$ respectively. The Langlands philosophy then {\it predicts} that there exists an isobaric automorphic representation $r(\pi) \boxtimes r'(\pi')$ of GL$(dd', \A_F)$ such that
$$
L(s, r(\pi) \boxtimes r'(\pi')) \, = \, L(s, \pi \times \pi'; r \otimes r').
\leqno(3.1)
$$
When it is known to exist, the map $\pi \to r(\pi)$ will be called a {\it functorial transfer} attached to $r$; some also call it a lifting.This is far from being known in this generality, but nevertheless, there have been some notable instances of progress wwhich we will make use of. Sometimes we do not know $r(\pi)$ exists, but still one has some good properties of thee relevant $L$-functions.

\medskip

A cuspidal automorphic representation $\pi$ is said to be {\bf dihedral} iff  it admits a self-twist by a (necessarily) quadratic character $\delta$, i.e., $\pi \simeq \pi \otimes \delta$. Equivalently, there is a quadratic extension $K/F$ and a character $\chi$ of $K$, such that $\pi$ is isomorphic to $I_K^F(\chi)$, the representation {\it automorphically induced} by $\chi$ from $K$ (to $F$). The passage from the second to the first definition is by taking $\delta$ to be the quadratic character of $F$ associated to $K/F$.

\medskip

We will need to use the following results:

\medskip

\noindent{\bf Theorem 3.2} ([Ra2000]) \, \it Let $\pi, \pi'$ be cuspidal automorphic representations of GL$(2, \A_F)$. Then there exists an isobaric automorphic representation $\pi \boxtimes \pi'$ of GL$(4, \A_F)$ such that 
$$
L(s, \pi \boxtimes \pi') \, = \, L(s, \pi \times \pi').
$$
Moreover, $\pi \boxtimes \pi'$ is cuspidal iff one of the following happens:
\begin{itemize}
\item[(i)] $\pi, \pi'$ are both non-dihedral {\it and} there is no character $\mu$ such that $\pi' \simeq \pi \otimes \mu$;
\item[(ii)] One of them, say $\pi'$, is dihedral, with $\pi' = I_K^F(\chi)$ for a character $\chi$ of a quadratic 
extension $K$, and
the base change $\pi_K$ is cuspidal and not
isomorphic to $\pi_K \otimes (\mu \circ \theta)\mu^{-1},$
where $\theta$ denotes the non-trivial automorphism of $K/F$.
\end{itemize}
\rm

\medskip

Note that in case (ii), $\pi$ may or may not be dihedral, and in the latter situation, $\pi \boxtimes \pi'$ is cuspidal.

\medskip

If $L(s) = \prod_v L_v(s)$ is an Euler product, and if $T$ is a finite set of places, we will write $L^T(s)$ for the incomplete Euler product $\prod_{v \notin T} L_v(s)$.

\medskip

\noindent{\bf Theorem 3.3} ([GJ79] for $n=2$, [PPS89] for $n=3$ and [BuG92] for general $n$) \, \it Let $\pi$ be a cuspidal automorphic representation of GL$(n, \A_F)$. Let $S$ be the union of the archimedean places of $F$ with the set of finite places where $\pi$ is ramified. Then $L^S(s, \pi; {\rm sym}^2)$ admits a meromorphic continuation and is holomorphic in in the real interval $(1/2,1)$.
\rm

\medskip

When $n=2$, there is even an isobaric automorphic representation sym$^2(\pi)$ of GL$(3, \A_F)$ such that 
$$
L(s, {\rm sym}^2(\pi)) \, = \, L(s, \pi; {\rm sym}^2) \quad \quad (n=2),
$$
and sym$^2(\pi)$ is cuspidal iff $\pi$ is non-dihedral. 
\rm

\medskip

We are stating here only the facts which we need. The reader is urged to read the articles quoted to get the full statements. The functional equation and the meromorphic continuuation of the symmetric square $L$-functions of GL$(n)/F$ can also be deduced from the Langalnds-Shahidi method.

\medskip

\noindent{\bf Theorem 3.5} ([KSh2000], [K2000], [KSh2001]) \, Let $\pi$ be a cuspidal automorphic representation of GL$(2, \A_F)$. Then for $j = 3,4$, there is an isobaric automorphic representation sym$^j(\pi)$ such that
$$
L(s, {\rm sym}^j(\pi)) \, = \, L(s, \pi; {\rm sym}^j) \quad \quad (j=3,4).
$$
Moreover, sym$^3(\pi)$ is cuspidal iff sym$^2(\pi)$ is cuspidal and does not admit a self-twist by a cubic character, while sym$^4(\pi)$ is cuspidal iff sym$^3(\pi)$ is cuspidal and does not admit a self-twist by a quadratic character.
\rm

\medskip

A cuspidal automorphic representation $\pi$ of GL$(2, \A_F)$ is said to be {\bf tetrahedral}, resp. {\bf octahedral}, iff sym$^2(\pi)$, resp. sym$^3(\pi)$, is cuspidal and admits a non-trivial self-twist by a cubic, resp. quadratic character.
We will say that $\pi$ is of {\bf solvable polydedral type} 
iff it is either dihedral or tetrahedral or octahedral.

\vskip 0.2in

\section{\bf Proof of Theorem A}

\bigskip

In this section we will say that a pair $(\pi, \pi')$ of cuspidal automorphic representations of GL$(2, \A_F)$ is of {\bf general type} iff we have:

\noindent{$(4.1)$}
\begin{itemize}
\item[(a)] Neither $\pi$ nor $\pi$ is dihedral; \, and
\item[(b)] $\pi'$ is not a twist of $\pi$ by a character.
\end{itemize}

\medskip

First we will deal with the special cases when (a) or (b) does not hold.

\medskip

Suppose (a) is satisfied, but not (b), i.e., there is a chharacter $\mu$ of (the idele class group of) $F$ such that 
$$
\pi' \, \simeq \, \pi \otimes \mu.
$$ 
Then we have the decomposition
$$
\pi \boxtimes \pi' \, \simeq \, ({\rm sym}^2(\pi) \otimes \mu) \boxplus \omega\mu,
\leqno(4.2)
$$
where $\pi \boxtimes \pi'$ denotes the isobaric automorphic representation of GL$(4, \A_F)$ associated to $(\pi, \pi')$ in [Ra2000], and $\omega$ is the central character of $\pi$. In terms of $L$-functions, we have
$$
L(s, \pi \times \pi') \, = \, L(s, {\rm sym}^2(\pi) \otimes \mu)L(s, \omega\mu).
\leqno(4.3)
$$
One knows that , since $\pi$ is non-dihdral, $L(s, $sym$^2(\pi) \otimes \mu)$ admits no Landau-Siegel zero. This was proved in the ground-breaking article [GHLL94] for $\mu =1$ and $\pi$ self-dual; the general case was taken care of by a combination of the arguments of [HRa95] and then [Ba97]. Besides, when $\omega\mu$ is not self-dual, i.e., not of order $\leq 2$, $L(s, \omega\mu)$ admits no Siegel zero (see for example [HRa95]). Finally, it is a well known classical fact that for any character $\chi$ of order $\leq 2$, $L(s, \chi)$ can have at most one Siegel zero. So, putting all this together, we see that 

\noindent{$(4.4)$}
\begin{itemize}
\item[($\alpha$)]  The  Landau-Siegel zeros of $L(s, \pi \times (\pi \otimes \mu))$ coincide with those of $L(s, \omega\mu)$, \, and
\item[($\beta$)]  This set is empty unless $\omega\mu$ is of order $\leq 2$, in which case there is at most one Landau-Siegel zero.
\end{itemize}
If $F$ is a Galois number field (over $\Q$) not containing any quadratic field, one knows by [Stk74] that the Dedekind zeta function of $F$ has no Landau-Siegel zero. So we may replace {\it order $\leq 2$} in $(\beta)$ by {\it order $2$} for such $F$. This gives the desired assertion in this case, and it also brings up case (i) of Theorem A.

\medskip

Next consider the case when $\pi$ is non-dihedral, but $\pi'$ is dihedral, associated to a chaaracter $\chi$ of a quadratic extension $K$ of $F$. We will write $\pi' = I_K^F(\chi)$ and say that it is automorphically induced from $K$ to $F$ by $\chi$. Then by the basic properties of base change ([AC89], [Ra2000]) we have
$$
\pi \boxtimes \pi' \, \simeq \, I_K^F(\pi_K \otimes \chi),
\leqno(4.5)
$$
where $\pi_K$ denotes the base change of $\pi$ to GL$(2)/K$, which is cuspidal because $\pi$  is non-dihedral. Thus by the inductive nature of $L$-functions, we get the following identity:
$$
L(s, \pi \times \pi') \, = \, L(s, \pi_K \otimes \chi),
\leqno(4.6)
$$
By [HRa95] we know that $L(s, \pi_K \otimes \chi)$ does not admit any Landau-Siegel zero, and this gives Theorem A in this case.

\medskip

Now suppose both $\pi$, $\pi'$ are both dihedral. Then $\pi$, resp. $\pi'$, is naturally attached to a dihedral representation $\sigma$, resp. $\sigma'$, of the global Weil group $W_F$. Say, $\sigma = $Ind$_K^F(\chi)$, for a character $\chi$ of the Weil group of a quadratic extenssion. (By abuse of notation, we are writing Ind$_K^F$ instead of Ind$_{W_K}^{W_F}$. Since $\boxtimes$ corresponds to the usual tensor product on the Weil group side (see [Ra2000]), we have 
$$
L(s, \pi \times \pi') \, = \, L(s, {\rm Ind}_K^F(\chi) \otimes \sigma').
\leqno(4.7)
$$
By Mackey,
$$
{\rm Ind}_K^F(\chi) \otimes \sigma' \, \simeq \, 
{\rm Ind}_{K}^{F}(\chi \otimes {\rm Res}_K^F(\sigma')),
\leqno(4.8)
$$
where Res$_K^F$ denotes the restriction functor taking representations of $W_F$ to ones of $W_F$.

Suppose $\sigma'$ is also {\it not} induced by a character  of $W_K$. Then
${\rm Res}_K^F(\sigma')$ is irreducible and the base change $\pi'_K$ is cuspidal, and since $L(s, \pi \times \pi')$ equals $L(s, \pi'_K \otimes \chi)$, it has no Landau-Siegel zero, thanks to [HRa95].

So we may assume that $\sigma'$ is also induced by a character $\chi'$ of $W_K$. Then
$$
{\rm Res}_K^F(\sigma') \, \simeq \, \chi' \oplus (\chi' \circ \theta),
$$
where $\theta$ denotes the non-trivial automrophism of $K/F$. 
Plugging this into (4.8) and making use of the inductive nature f $L$-functions, we get
$$
L(s, \pi \times \pi') \, = \, L(s, \chi\chi')L(s, \chi(\chi' \circ \theta)).
\leqno(4.9)
$$
So there is no Landau-Siegel zero unless $\chi\chi'$ or $\chi(\chi' \circ \theta)$ is of order $\leq 2$, which we will asssume to be the case from now on. We have yet to show that there is at most one Landau-Siegel zero , which is true (see the remarks above) if only one of them has order $\leq 1$. 
Suppose they are both of order $\leq 2$. If one of them, say $\chi\chi'$ is trivial, then 
$$
L(s, \pi \times \pi') \, = \, \zeta_F(s)L(s, \nu),
\leqno(4.10)
$$  
where $\nu = \frac{\chi' \circ \theta}{\chi'}$. Note that since $\sigma'$ is irreduucible, $\chi'$ is not equal to $\chi' \circ \theta$.Then $\nu$ must be a quadratic character, and the right hand side of (4.10) evidently defines a non-negative Dirichlet series with a pole of order $1$ at $s=1$. So by  Lemma 1.7, $L(s, \pi \times \pi')$ can have at most one Landau-Siegel zero.

It is left to consider when $\mu: = \chi\chi'$ and $\nu'$ are both quadratic characters. The argument here is well known, and we give it nly for thre sake of completeness. Notes that the Dirichlet series defined by
$$
L(s): \, = \, \zeta_F(s)L(s, \mu)L(s, \nu)L(s, \mu\nu),
\leqno(4.11)
$$
has non-negative coefficieents, meromorphic continuation and a functional equation, with no pole except at $s=1$, where the pole is simple; $L(s)$ is the Dedekind zeta function of the biquadratic extension of $F$ obtained as the compositum of the quadratic extensions cut out by $\mu$ and $\nu$. Thus by applying Lemma 1.7 again, we see that $L(s)$, and hence its divisor $L(s, \pi \times \pi')$ (see (4.9), has at most one Landau-Siegel zero.

This finishes the proof of Theorem A when $\pi, \pi'$ are both dihedral, bringing up case (ii) when they are both defined by chracters of the same quadratic extension $K$.

\medskip

So we may, and we will, asssume from here on that bothe (a) and (b) of (4.1) are satisfied. Now Theorem A will be proved if we establish the following theorem, which gives a stronger statement.

\medskip

\noindent{\bf Theorem 4.12} \, \it
	Let $\pi$ and $\pi'$ are unitary cuspidal automorphic representations of GL$(2, \A_F)$,and assume that the pair $(\pi, \pi')$ is of general type. Then,

	\textnormal{(a)}: There is an effective absolute constant $c \gneqq 0$
	such that $L(s, \pi \times \pi')$ has no zero in the interval
	$(\, 1 - c / {\log M}, 1\,)$.

	\textnormal{(b)}: Additionally, if $\pi$ and $\bar{\pi}'$ is not twist
	equivalent by a product of a quadratic character and ${|\,|}^{\fI t}$,
	then there exists an absolute effective constant $c_{2} \gneqq 0$
	such that $L(s, \pi \times \pi')$ has no zero in the region

	\[
		\Set{s = \sigma + \fI t \mid \sigma \le 1 - {(c_{2} \fL_{t})}^{-1}}
	\]

	where

	\[
		\fL_{t} = \log{[N(\pi \times \pi') \, D_{F}^{4} \,
			 {(2 + |t| + \Lambda)}^{4 N}]},
	\]
with $N = [F : \mQ]$ and $\Lambda$ denoting the maximum of the infinite types
of $\pi$ and $\pi'$.
\rm

\medskip

See 1.4 for the definition of the infinite parameter $\Lambda$. Such a result was
a working hypothesis in the work of Moreno ([Mo]) on a n effective vesion of the strong multiplicity one theorem for GL$(2)$.

\noindent
\emph{Proof. of Theorem 4.12}

(a) \, Put
$$
\fL \, = \, \fL_0.
$$
Then by definition of $\Lambda$, 
$$
\log M \, = \, \fL.
\leqno(4.13)
$$

\medskip 

Let $\omega$ and $\omega'$ be the central characters of $\pi$ and $\pi'$
respectively.

Since $\pi, \pi'$ is nondihedral, sym$^2(\pi)$ and ${\rm sym}^{2} (\pi')$
are cuspidal. Also, $(\pi, \pi')$ being of general type implies 
(cf. [Ra2000]) that their Rankin-Selberg product $\pi \boxtimes \pi'$ of GL$(4, \A_F)$ is cuspidal.

Consider the following isobaric automorphic representation
$$
	\Pi = 1 \boxplus (\pi \boxtimes \pi') \boxplus
	({\rm sym}^{2} (\bar{\pi}) \otimes \omega)
\leqno(4.14)
$$
Write, as usual
$$
\Pi \, = \, \Pi_\infty \otimes \Pi_f.
$$
Note that $\Pi$ is unitary and so its contragredient $\Pi^\vee$ identifies with its complex conjugate $\bar \Pi$.
By the bi-additivity of the Rankin-Selberg process, we have the factorization

\noindent{$(4.15)$}
\begin{align}
L (s, \Pi_f \times \bar{\Pi}_f) = &\zeta_{F} (s) L (s, \pi_f \times \pi'_f)
		L (s, \bar{\pi}_f \times \bar{\pi}'_f) 
	\notag \\
	&
	L (s, {\rm sym}^{2} (\pi_f) \otimes {\omega}^{-1})
	\notag \\	
	&
	L (s, {\rm sym}^{2} (\bar{\pi_f}) \otimes \omega)
	L (s, (\pi_f \boxtimes \pi'_f) \times (\bar{\pi}_f \boxtimes \bar{\pi}'_f))
	\notag \\
	&
	L (s, {\rm sym}^{2} (\pi_f) \times {\rm sym}^{2} (\bar{\pi}_f))
	L (s, (\pi_f \boxtimes \pi'_f) \times {\rm sym}^{2} (\pi_f) \otimes {\omega}^{-1})
	\notag \\
	&
	L (s, (\bar{\pi}_f \boxtimes \bar{\pi}'_f) \times {\rm sym}^{2} (\bar{\pi}_f) \otimes \omega).
	\notag
\end{align}
By abuse of notation, we are writing $\omega$ instead of $\omega_f$, which should not cause any confusion.

It is well known that $\zeta_F(s)$ has a simple pole at $s=1$, and since sym$^2(\pi)$ and $\pi \boxtimes \pi'$ are cuspidal, $L (s, (\pi_f \boxtimes \pi'_f) \times (\bar{\pi_f} \boxtimes \bar{\pi}'_f))$ and $L(s, {\rm sym}^{2} (\pi_f) \times {\rm sym}^{2} (\bar{\pi}_f))$ have simple poles at $s=1$ as well. Moreover, the remaining factors are entire with no zero at $s=1$ (see the discussion following (2.3)). Thus
$$
{\rm ord}_{s=1} L(s, \Pi_f \times \bar \Pi_f) \, = \, 3.
\leqno(4.16)
$$ 

By Lemma 2.5, the Dirichlet series defined by $L (s, \Pi_f \times \bar{\Pi_f})$
has nonnegative coefficients.

Put
$$
	L_{1} (s) = L (s, \pi_f \times \pi'_f)  L (s, \bar{\pi}_f \times \bar{\pi}'_f)
\leqno(4.17)
$$
Since the real zeros of  $L (s, \pi_f \times \pi'_f)$ and $ L (s, \bar{\pi}_f \times \bar{\pi}'_f)$ are the same, we get for any $\beta \in (0,1)$,
$$
{\rm ord}_{s=\beta} L_1(s) \, = \, 2 \, {\rm ord}_{s=\beta} L(s, \pi_f \times \pi'_f).
\leqno(4.18)
$$

Next observe that at any place $v$, if $\sigma_v$ (resp. $\sigma'_v$) denotes the $2$-dimensional representation of $W'_{F_v}$ (resp. $W_{F_v}$) attached to $\pi_v$ for $v$ finiite (resp. $v$ archimedean) by the local Langlands correspondence, we have
$$
\sigma_v \otimes {\rm sym}^3(\sigma_v) \, \simeq \, (\sigma_v \otimes \omega_v) \oplus {\rm sym}^3(\sigma_v),
$$
which implies the decomposition 

\noindent{$(4.19)$}
\begin{align}
(\sigma_v \otimes \sigma'_v) \otimes {\rm sym}^2(\sigma_v) \otimes \omega_v^{-1}  &\simeq
(\sigma_v \otimes {\rm sym}^2(\sigma_v) \otimes \omega_v^{-1}) \otimes \sigma_v'
\notag \\
&\simeq (\sigma_v \otimes \sigma'_v) \oplus ({\rm sym}^3(\sigma_v) \otimes \omega_v^{-1} \otimes \sigma'_v).
\end{align}

This gives, by the definition of automorphic $L$-functions in section 1, the following identity of $L$-functions:
$$
L (s, (\pi_f \boxtimes \pi'_f) \times {\rm sym}^{2} (\pi_f) \otimes {\omega}^{-1})
\, = \, L (s, \pi_f \times \pi'_f) L (s, A^3(\pi_f) \times \pi'_f)
\leqno(4.20)
$$
where, following [KSh2001], we have set
$$
A^3(\pi) : \, = \, {\rm sym}^3(\pi)	\otimes {\omega}^{-1}.
$$
We need

\medskip

\noindent{\bf Lemma 4.21} \, \it Since $(\pi, \pi')$ is of general type, $L(s, A^3(\pi_f) \times \pi_f')$ and $L(s, A^3(\bar \pi_f) \times \bar \pi'_f)$ are entire.
\rm

\medskip

{\it Proof of Lemma} \, Existence of a pole for one of them, say at $s=s_0$, will imply a pole for the other at $s = \overline s_0$; hence it suffices to prove that $L(s, A^3(\pi_f) \times \pi'_f)$ is entire. Since the local factors at the archimedean places do not vanish, it is enough to show that the full $L$-function $L(s, A^3(\pi) \times \pi')$ is entire. Since $(\pi, \pi')$ is of general type, $\pi, \pi'$ are non-dihedral and not twists of each ther. If $\pi$ is not tetrahedral (see section 3 for definition), then by [KSh2000], sym$^3(\pi)$ is cuspidal. The asserton of Lemma is clear in thhis case by the standard results on the Rankin-Selberg $L$-functions (see section 2). So we may, and we will, assume that $\pi$ is tetrahedral. Then sym$^2(\pi)$ is isomorphic to sym$^2(\pi) \otimes \nu$ for some cubic character $\nu$, and by Theorem 2.2 of [KSh2001], $A^3(\pi)$ is isomorphic to $(\pi \otimes \nu) \boxplus (\pi \otimes \nu^2)$. Then $L(s, A^3(\pi) \times \pi')$ factors as $L(s, (\pi \otimes \nu) \times \pi')L(s, (\pi \otimes \nu^2) \times \pi')$,
which is entire by the Rankin-Selberg theory, because $\pi'$ is not a twist of $\bar \pi \simeq \pi \otimes \omega^{-1}$.

\qed

\medskip

Put

\noindent{$(4.22)$}

\begin{align}
L_{2} (s)= &\zeta_{F} (s) L (s, {\rm sym}^{2} (\pi_f) \otimes {\omega}^{-1})	
L (s, {\rm sym}^{2} (\bar{\pi_f}) \otimes \omega)
	\notag \\
	&L (s, (\pi_f \boxtimes \pi'_f) \times (\bar{\pi}_f \boxtimes \bar{\pi}'_f))
	L (s, {\rm sym}^{2} (\pi_f) \times {\rm sym}^{2} (\bar{\pi}_f))
	\notag \\
	&L (s, A_3(\pi_f) \times \pi'_f) 
	L (s, \bar A_3(\pi)_f \times \bar{\pi}_f')
	\notag
\end{align}

Then
$$
L (s, \Pi \times \bar{\Pi}) \, = \,	L_{1}^{2} (s) L_{2} (s). 
\leqno(4.23)
$$

Applying Lemma 4.21, and using the cuspidality of sym$^2(\pi)$ and $\pi \boxtimes \pi'$, we get the
following

\medskip

\noindent{\bf Lemma 4.24} \, \it $L_2(s)$ is entire.
\rm

\medskip
 
Combining this lemma with (4.16), (4.23) and (2.x), and using Lemma 1.7, we get the
existence of a positive, effective constant $c$ such that 
$$
2 {\rm ord}_{s=\beta} L_{1}(s) \, \leq \, 3 \quad if \quad \beta \in (1-c/\log M, 1).
\leqno(4.25)
$$ 
In view of (4.18), if $L (s, \pi \times \pi')$
has a Landau-Siegel zero $\beta$ (relative to $c$), 
then $\beta$ will be a zero of
$L_{1}^{2} (s)$ of multiplicity $4$, leading to a contradiction.

We have now proved part (a) of Theorem 4.12, and hence Theorem A.

\bigskip

(b) \, First note that under the condition of (2),
$L(s, \pi \times \pi' \tensor | |^{\fI t})$ has no Landau-Siegel zero.
Moreover, the maximum $\Lambda_{t}$ of infinite types $\Lambda(\pi)$
and $\Lambda(\pi' \tensor ||^{\fI t})$ are no more than $|t| + \Lambda$.
Thus $L(s, \pi \tensor ||^{\fI t})$ has no zero in the interval
$$
	1 - \frac{1}{(c_{2} L_{t})} < \sigma < 1
\leqno(4.26)
$$
Since we have
$$
L(\sigma + \fI t, \pi \times \pi') \, = \, L(\sigma, \pi \times \pi' \tensor ||^{\fI t}),
$$
the assertion of (b) now follows. 

\qedsymbol 

\vskip 0.2in

\section{\bf Proof of Theorem B}

\bigskip

Let $\pi$ be a 
cuspidal automorphic representation of GL$(2, \A_F)$ of central character $\omega$. First we will dispose of the {\it solvable polyhedral} cases, where we will not need to assume that $\pi$ is self-dual. 

Suppose $\pi$ is {\it dihedral}, i.e., of the form $I_K^F(\chi)$ for a character $\chi$ (of the idele classes) of a quadratic extension $K$ of $F$, with  $\theta$ denoting non-trivial automorphism of $K/F$.  
Let $\chi_o$ denote the restriction of $\chi$ to $F$. Note that
$$
\chi\chi^\theta \, = \, \chi_0 \circ N_{K/F},
\leqno(5.1)
$$
where $N_{K/F}$ denotes the norm from $K$ to $F$. ($\chi_0 \circ N_{K/F}$ is the base change $(\chi_0)_K$ of $\chi_0$ to $K$.) In particular,
$$
I_K^F(\chi\chi^\theta) \, \simeq \, \chi_0 \boxplus \chi_0\delta,
\leqno(5.2)
$$
where $\delta$ denotes the quadratic character of $F$ attached to $K/F$ by class field theory. For any pair $(\lambda, \xi)$ of characters of $K$, one has (cf. [Ra2000])
$$
I_K^F(\lambda) \boxtimes I_K^F(\xi)  \, \simeq \, I_K^F(\lambda\xi) \boxplus I_K^F(\lambda\xi^\theta).
\leqno(5.3)
$$
Putting $\lambda = \xi = \chi$ in (5.3), and using (5.1), (5.2), we get
$$
\pi \boxtimes \pi \, \simeq \, I_K^F(\chi^2) \boxplus \chi_0 \boxplus \chi_0\delta.
$$
Since $\pi \boxtimes \pi$ is the isobaric sum ($\boxplus$) of sym$^2(\pi)$ with  $\omega$, which is $\chi_0\delta$ (as it corresponds to the determinant of the representation Ind$_K^F(\chi)$ of $W_K$), we get
 $$
{\rm sym}^2(\pi)  \, \simeq \, I_K^F(\chi^2) \boxplus \chi_0,
\leqno(5.4)
$$
Putting $\lambda = \xi = \chi^2$ in (5.3), using (5.1), (5.2), (5.4), and the inductive nature of $L$-functions, we get the following identity of $L$-functions:
$$
L(s, {\rm sym}^2(\pi) \times {\rm sym}^2(\pi)) \, = \, L(s, \chi^4)L(s, \chi_0^2)^2 L(s, \chi_0^2\delta) L(s, \chi^3\chi^\theta)^2.
\leqno(5.5)
$$
It is an abelian $L$-function, and the problem of Landau-Siegel zeros here is the classical one, and there is no such zero unless one (or more) of the characters appearing on the right of (5.5) is of order $\leq 2$. When $\omega = 1$, $\chi_0$ is $\delta$, and since $\delta^2 =1 = \delta \circ N_{K/F}$, we obtain
$$
L(s, {\rm sym}^2(\pi) \times {\rm sym}^2(\pi)) \, = \, L(s, \chi^4)\zeta_F(s)^2L(s, \delta)L(s, \chi)^2.
\leqno(5.6)
$$

Next let $\pi$ be {\it tetrahedral}, in which case sym$^2(\pi)$ is cuspidal and admits a self-twist by  a non-trivial cubic character $\mu$.
In other words, there is a cyclic extension $M/F$ of degree $3$ cut out by $\mu$, with non-trivial automorphism $\alpha$, and a character $\lambda$ of $M$, not fixed by $\alpha$, such that
$$
{\rm sym}^2(\pi) \, \simeq \, I_M^F(\lambda).
\leqno(5.7)
$$
Since by Mackey,
$$
{\rm Ind}_M^F(\lambda)^{\otimes 2} \, \simeq \, {\rm Ind}_M^F(\lambda^2) \oplus {\rm Ind}_M^F(\lambda\lambda^\alpha) \oplus {\rm Ind}_M^F(\lambda\lambda^{\alpha^2})
$$
we get
$$
L(s,  {\rm sym}^2(\pi) \times {\rm sym}^2(\pi)) \, = \, L(s, \lambda^2) L(s, \lambda\lambda^\alpha)
L(s, \lambda\lambda^{\alpha^2}).
\leqno(5.8)
$$
Again it is an abelian $L$-function, and there is nothing more to prove.

Now let $\pi$ be {\it octahedral}. Then by definition, sym$^j(\pi)$ is cuspidal for $j \leq 3$ and moreover,
$$
{\rm sym}^3(\pi) \, \simeq \, {\rm sym}^3(\pi) \otimes \eta,
\leqno(5.9)
$$
for a quadratic character $\eta$. Equivalently, there is a quadratic extension $E/F$ (attached to $\eta$) such that the base change $\pi_E$ is tetrahedral, ie., there exists a cubic character $\nu$ of $E$ such that
$$
{\rm sym}^2(\pi_E) \, \simeq \, {\rm sym}^2(\pi_E) \otimes \nu.
\leqno(5.10)
$$
Now we appeal to the evident identity
$$
L(s,  {\rm sym}^2(\pi) \times {\rm sym}^2(\pi)) \, = \, L(s, {\rm sym}^4(\pi))L(s, {\rm sym}^2(\pi)\otimes \omega) \zeta_F(s).
\leqno(5.11)
$$
Then it suffices to prove that the set of Landau-Siegel zeros of $L(s, {\rm sym}^4(\pi))$ is the same as that of the maximal abelian $L$-function dividing it. To this end we note that by Theorem 3.3.7, part (3), of [KSh2001], 
$$
{\rm sym}^4(\pi) \, \simeq \, I_E^F(\nu^2) \otimes \omega^2 \boxplus {\rm sym}^2(\pi) \otimes \omega\eta,
\leqno(5.12)
$$
so that
$$
L(s, {\rm sym}^4(\pi)) \, = \, L(s, \nu^2(\omega \circ N_{E/F})^2)L(s,  {\rm sym}^2(\pi) \otimes \omega\eta).
\leqno(5.13)
$$
Recall from section 1 that since $\pi$ is non-diehdral, $L(s, {\rm sym}^2(\pi) \otimes \beta)$ has no Landau-Siegel zeero for any character $\beta$. So we are done in this case as well.

\medskip

So from now on we may, and we will, assume that $\pi$ is {\it not of solvable polyhedral type}. In view of the identity (5.11), the derivation of Theorem B will be complete once we prove the following result on the symmetric $4$-th power $L$-function of $\pi$, which may be of independent interest.

\medskip

\noindent{\bf Theorem B$^\prime$} \, \it
Let $\pi$ be a cusspidal automorphic represenation of GL$(2, \A_F)$ with trivial central  character, which is not of solvable polyhedral type. Then $L(s, {\rm sym}^{4} (\pi))$ admits no Landau-Siegel zero, 
	More explicitly, there exists a postive, effective constant $C$ such that it has no zero in the real interval
	$(\,1 - C {\fL}^{-1} \,)$
	for some constant $C$, where
	\[
		\fL = \log [N(\pi) D_{F}^{2} {(2 + \Lambda)}^{2 N}]
	\]
	where $N = [F : \mQ]$, and $\Lambda$  the infinite type of $\pi$.
\rm

\medskip

\noindent{\bf Corollary 5.14}  \, \it
	Under the hypotheses of Theorem B$^\prime$, the Landau-Siegel zeros of
	$L(s, {\rm sym}^{2} (\pi) \times {\rm sym}^{2} (\pi))$, if any, comes from one of
 $\zeta_F(s)$. If $F$ is a Galois extension of $\Q$ not containing any quadratic field, there is no Landau-Siegel zero at all.
\rm

\medskip
  
\noindent
\emph{Proof of Theorem B$^\prime$}

First we note that ${\rm sym}^{4} (\pi)$ is cuspidal as $\pi$ is not of solvable polyhedral type ([KSh2001]).
Also, ${\rm sym}^4(\pi)$ is self-dual as $\pi$ is.

Put 
$$
\Pi = 1 \boxplus {\rm sym}^{2} (\pi)  \boxplus {\rm sym}^{4} (\pi),
\leqno(5.15)
$$
which is a self-dual isobaric automorphic representation of GL$(9, \A_F)$.
Since it is unitary, it is also self-conjugate.

A formal caculation gives the identities
$$
L(s, {\rm sym}^{4}(\pi) \times {\rm sym}^{2} (\pi)) \, = \,
L(s, {\rm sym}^{2} (\pi))L(s, \pi; {\rm sym}^6),
\leqno(5.16)
$$
and

\noindent{$(5.17)$}
\begin{align}
L(s, \Pi \times \Pi) =
&\zeta_{F} (s) L(s, {\rm sym}^{2} (\pi) \times {\rm sym}^{2} (\pi))
L(s, {\rm sym}^{4} (\pi) \times {\rm sym}^{4} (\pi))
\notag \\
&\phantom{=} L(s, {\rm sym}^{2} (\pi))^2
L(s, {\rm sym}^{4} (\pi))^{2}
L(s, {\rm sym}^{4}(\pi) \times {\rm sym}^{2} (\pi))^2
	\notag \\
	&= \zeta_{F} (s) L(s, {\rm sym}^{2} (\pi) \times {\rm sym}^{2} (\pi))
	L(s, {\rm sym}^{4} (\pi) \times {\rm sym}^{4} (\pi))
	\notag \\
	&\phantom{=} 
	 L(s, {\rm sym}^{2} (\pi))^4
	L(s, {\rm sym}^{4} (\pi))^{4} L(s, \pi; {\rm sym}^6)^2.
	\notag
\end{align}

By Lemma 2.5, the Dirichlet series defined by $L(s, \Pi_f \times \Pi_f)$ has non-negativ coefficieents and moreover, the cusspidality of sym$^j(\pi)$ for $j = 2,4$ implies that
$$
-{\rm ord}_{s=1} L(s, \Pi_f \times \Pi_f) \, = \, 3.
\leqno(5.18)
$$

Put 
$$
L_{1} (s) \, = \, L(s, {\rm sym}^{4} (\pi_f))^4
\leqno(5.19)
$$ 
and define $L_{2} (s)$ by the equation
$$
L(s, \Pi_f \times \Pi_f) \, = \, L_1(s)L_2(s).
\leqno(5.20)
$$

\medskip

\noindent{\bf Proposition 5.21} \, \it  $L_{2} (s)$ is holomorphic
in the interval $(\,1/2, 1\,)$.
\rm

\medskip

{\it Proof of Proposition} \,  Since we have

\noindent{$(5.22)$}
\begin{align}
	L_{2} (s) &= \zeta_{F} (s) L(s, {\rm sym}^{2} (\pi_f) \times {\rm sym}^{2} (\pi_f))
	L(s, {\rm sym}^{4} (\pi_f) \times {\rm sym}^{4} (\pi_f))
	\notag \\
	&\phantom{=} L(s, {\rm sym}^{2} (\pi_f))^4
	L(s, \pi_f; {\rm sym}^{6})^2,
	\notag
\end{align}
and since all the factors other than the square of the  symmetric $6$-th power $L$-function are, owing to the cuspidality of sym$^j(\pi)$ for $j = 2,4$, holomorphic in $(0,1)$, it suffices to show the same for $L(s, \pi_f; {\rm sym}^6)$. But this we cannot do, given the current state of what one knows.

But we are thankfully rescued by the following identity
$$
	L(s, {\rm sym}^{3} (\pi_f); {\rm sym}^{2})
	= L(s, {\rm sym}^{2} (\pi_f) ) L(s, \pi_f; {\rm sym}^{6}).
\leqno(5.23)
$$
Consequently, we have

\noindent{$(5.24)$}
\begin{align}
	L_{2} (s) &= \zeta_{F} (s) L(s, {\rm sym}^{2} (\pi_f) \times {\rm sym}^{2} (\pi_f))
	L(s, {\rm sym}^{4} (\pi_f) \times {\rm sym}^{4} (\pi_f))
	\notag \\
	&\phantom{=} L(s, {\rm sym}^{2} (\pi_f))^2
	L(s, {\rm sym}^3(\pi_f); {\rm sym}^2)^2,
	\notag
\end{align}
and Proposition 5.21 will follow from

\medskip

\noindent{\bf Lemma 5.25} \, \it $L(s, {\rm sym}^3(\pi_f); {\rm sym}^2)$ is holomorphic in $(1/2, 1)$.
\rm

\medskip

{\it Proof of Lemma 5.25} \, 
Let $S$ be the union of the archimedean places of $F$ with the (finite) set of fiinite places $v$ where ${\rm sym}^4(\pi)$ is unramified. It will be proved in section 7 (see Lemmas 7.9 and 7.4) that at any $v$ in $S$, $L(s, \pi_v; {\rm sym}^{2j})$ is holomorphic in $(1/2, 1)$ for $j \leq 4$. 

Thus it suffices to prove that the incomplete $L$-function $L^S(s, {\rm sym}^3(\pi); {\rm sym}^2)$, defined in a right half plane by $\prod_{v \notin S} L(s, {\rm sym}^3(\pi_v); {\rm sym}^2)$, is holomorphic in $(1/2, 1)$. But since sym$^3(\pi)$ is a cuspidal automorphic representation of GL$(4, \A_F)$, this is a consequence (see Theorem 3.3) of the main result of [BuG92].

\qed

\medskip

{\it Proof of Theorem B$^\prime$} (contd.) \,
Apply Lemma 1.7 and Proposition 2.9
to the positive Dirichlet series $L(s, \Pi_f \times \Pi_f)$, which has a pole at $s=1$ of order $3$. Since $L_2(s)$ is holomorphic in $(1/2,1)$, there is an effective constant $c > 0$ such that the number of real zeros of $L_1(s)$ 
in $(1 - \frac{c}{\log M(\Pi \times \Pi)}, 1)$ is bounded above by $3$. But $L_1(s)$ is the fourth power of $L(s, \pi_f; {\rm sym}^4)$, and so $L(s, \pi_f; {\rm sym}^4)$ can have no zero in this interval. Also, by (2.8), 
$$
M(\pi; {\rm sym}^4) \, \asymp \, M(\Pi \times \Pi),
\leqno(5.26)
$$
where the implied constants are effective.

Now we have proved Theorem B$^\prime$, and hence Theorem B.

\qedsymbol

\noindent{\bf Remark 5.27} \, In Theorem B, we asssumed that $\pi$ is self-dual. To treat the general case with these arguments one needs the following hypothesis for $r=4$.

\medskip

\noindent{\bf Hypothesis 5.28} \, \it
Let $\pi$ be a unitary cuspidal representaion of $GL(r, \A_F)$, and $\chi$
a non-trivial quadratic character of $F$,
then $L(s, \pi; {\rm sym}^{2} \otimes \chi)$ is holomorphic in
$(\,t, 1\,)$ for a fixed real number $t < 1$.
\rm

\medskip

For $r=2$, of course, there is nothing to do as sym$^2(\pi)$ is automorphic ([GJ79]). For $r=3$, this was established W.~Banks in [Ba96], thus proving a hypothesis of [HRa95] enabling the completion of the proof of the lack of Landau-Siegel zeros for cusp forms on GL$(3)/F$.

\vskip 0.2in

\section{\bf Proof of Corollary C}

\medskip

Here $g$ is a Maass form on the upper half plane relative to SL$(2, \Z)$, with Laplacian eigenvalue $\lambda$ and Hecke eigenvalues $a_p$. If $\pi$ is the cuspidal automorphic representation of GL$(2, \A)$, $\A = \A_\Q$, generated by $g$ 
(see [Ge75]), we may consider the Gelbart-Jacquet lift sym$^2(\pi)$, which is an isobaric automorphic representation
of GL$(3, \A)$. Since $g$ has level $1$, it is not dihedral, and so sym$^2(\pi)$ is cuspidal. Moreover, since sym$^2(\pi_p)$ is, for any prime $p$, {\it unramified} because $\pi_p$ is, which means that sym$^2(\pi_p)$ is {\it spherical} at $p$, i.e., it admits a non-zero vector fixed by the maximal compact subgroup $K_p : = $GL$(3, \Z_p)$. It is also spherical at infinity, i.e., has a non-zero fixed vector under  the orthogonal group $K_\infty: = {\rm O}(3)$. Moreover, the center $Z(\A)$ acts trivially, and the archimedean component sym$^2(\pi_\infty)$ consists of eigenfunctions for the center $\mathfrak z$ of the enveloping algebra of Lie$({\rm GL}(3, \R))$.  In sum, sym$^2(\pi)$ is a subrepresentation of 
$$
V: \, = \, L^2(Z(\A){\rm GL}(3, \Q) \backslash {\rm GL}(3, \A)),
\leqno(6.1)
$$
admitting a spherical vector, i.e., a (non-zero) smooth function $\phi$ invariant under $K: = \prod_v K_v$, where $v$ runs over
the places $\{\infty, 2,3,5,7, \ldots, p, \ldots\}$. GL$(3, \A)$ acts on $V$ by right translation and leaves invariant the natural sclar product $\langle \, .\, ,\, . \, \rangle$ given, for all $\xi_1, \xi_2 \in V$, by 
$$
\langle \, \xi_1 \, , \, \xi_2 \, \rangle \, = \, \int_{Z(\A){\rm GL}(3, \Q) \backslash {\rm GL}(3, \A)} \xi_1(x){\overline \xi_2}(x) dx,
$$
where $dx$ is the quotient measure defined by the Haar measures on ${\rm GL}(3, \A)$, $Z(\A)$ and ${\rm GL}(3, \Q)$, chosen as follows. On the additive group $\A$, take the measure to be the product measure $\prod_v dy_v$, where $dy_\infty$ is the Lebesgue measure on $\Q_\infty = \R$, and for each prime $p$, $dy_p$ is normalized to give measure $1$ to $\Z_p$. Take the measure $d^\ast y = dy/|y|$ on $\A^\ast$, where $|y| = \prod_v |y_v|$ the natural absolute value, namely the one given by $|y_\infty| = {\rm sgn}(y_\infty)y_\infty$ and $|y_p| = p^{-v_p(y_p)}$. Since the center $Z$ is isomorphic to the multiplicative group, this defines a Haar measure $dz = \prod_v dz_v$ on $Z(\A)$. On GL$(3, \A)$ take the product measure $\prod_v dx_v$, where each $dx_v$ is given, by using the Iwasawa decomposition GL$(3, \Q_v) = Z_vT_vN_vK_v$, as $dz_vdt_vdn_vdk_v$. Here $T_v$ denotes the subgroup of diagonal matrices of determinant $1$, with $dt_v$ being the transfer of the measure $d^\ast t_v$ via the isomorphism $T_v \simeq F_v^\ast$, $N_v$ the unipotent upper triangular group with measure $dn_v$ being the transfer of $dt_v$ via the isomorphism of $N_v$ with the additive group $F_v$, and $dk_v$ the Haar measure on $K_v$ noralized to give total volume $1$. 

\medskip
 
The representation sym$^2(\pi)$ is a unitary summand. Since sym$^2(\pi)$ is irreducible, such a $\phi$ will generate the whole space by taking linear combinations of its translates and closure. Note that $\phi$ is the pull back to GL$(3, \A)$ of a  function $\phi_0$, which is real analytic by virtue of being a $\mathfrak z$-eigenfunction, on the $5$-dimensional (real) orbifold
$$
M : \, = \,  Z(\A){\rm GL}(3, \Q) \backslash {\rm GL}(3, \A)/K \, = \, SL(3, \Z)\backslash {\rm SO}(3).
\leqno(6.2)
$$
Since $\phi$ and $\phi_0$ determine each other, we will by abuse of notation use the same symbol $\phi$ to denote both of them. 

\medskip

The spherical function $\phi$, sometimes called a {\it new vector}, is unique only up to multilication by a scalar. It is important for us to normalize it. There are two natural ways to do it. The first way, called the {\it arithmetic normalization}, is to make the Fourier coefficient $a(1,1)$ (see below) equal $1$ (as for newforms on the upper half plane). The second way, which is what we will pursue here, is called the {\it spectral normalization}, and normalizes the scalar product $\langle, \rangle$ of $\phi$ with itself to be essentially $1$.When so normalized, we will use the symbol sym$^2(f)$ for $\phi$. We will appeal to the Fourier expansions in terms of the Whittaker functions to do it. We begin with the general setup.

\medskip

\noindent \emph{Definition 6.3} Let $\phi$ be automorphic form on GL$(n)/\Q$ generating
a unitary, spherical, cuspidal automorphic representation $\Pi$.
Say that $\phi$ is \emph{normalized} if we have:

$$
	\phi (g) \, = \, \sum_{\gamma \in U(n - 1, \Q) \backslash {\rm GL}(n - 1, \Q)}
	W_{\Pi} \left(
	\begin{pmatrix}
		\gamma & \phantom{0} \\
		\phantom{0} & 1 
	\end{pmatrix}
	g \right),
\leqno(6.4)
$$
where $U(n - 1, \Q)$ denotes the subgroup of $GL(n - 1, \Q)$ consisting of
upper triangular, unipotent matrices, $W_{\Pi} = \prod_{v} W_{\Pi, v}$ the global Whittaker function whose local
components are defined below. (Again, $\Pi$ spherical means that $\Pi_v$ admits, at every place $v$, a non-zero vector invariant under the maximal compact (mod center) subgroup $K_v$, which is 
GL$(n, \mathfrak \Z_p)$ when $v$ is $v_p$ for a prime $p$.)

\medskip

\noindent{$(6.5)$}
\begin{itemize}
\item{$W_{\Pi, p}$ is, for any prime $p$, the unique $K_{p}$-invariant
function corresponding to $\Pi_{p}$ normalized so that $W_{\Pi, p} (e) = 1$.} 
\item{At the archimedean place,
$$
W_{\Pi, \infty} = C(\Pi)^{-\frac{1}{2}} W_{n, a},
$$
where $W_{n, a}$ be the normalized spherical function of infinite type $a$ on $GL(n, \R)$
in the sense of Stade [St2001], and 
$$
C(\Pi) = L(1, \Pi_{\infty} \times \Pi_{\infty}).
$$}
\end{itemize}

Denote the function so normalized in the space of $\Pi$ by the symbol $\phi(\Pi)$.

\medskip

Now let us get back to our Maass form $g$ for SL$(2, \Z)$, with  associated spherical cuspidal 
representation $\pi$, resp. 
${\rm sym}^{2} (\pi)$, of GL$(2, \A)$, resp. GL$(3, \A)$. We set
$$
{\rm sym}^2(g) \, = \, \phi({\rm sym}^2(\pi)).
\leqno(6.6)
$$
Since $g$ has level $1$, one knows that $\lambda > 1/4$ (in fact $> 50$, though we do not need it), so that if we write
$$
\lambda \, = \, \frac{1-t^2}{4}, 
$$
then $t$ is a non-zero real number; so we may choose $t$ to be positive. We have
$$
L(s, \pi_\infty) \, = \, \Gamma_\R(s+it)\Gamma_\R(s-it).
$$
Consequently,
$$
L(s, {\rm sym}^2(\pi_\infty)) \, = \, \Gamma_\R(s+2it)\Gamma_\R(s)\Gamma_\R(s-2it),
$$
and
$$
L(s, {\rm sym}^2(\pi_\infty) \times {\rm sym}^2(\pi_\infty)) \, = \, 
\Gamma_\R(s+4it)\Gamma_\R(s+2it)^2\Gamma_\R(s)^3\Gamma_\R(s-2it)^2\Gamma_\R(s-4it).
$$
Then, since $\Gamma(1-ait)$ is the complex conjugate of $\Gamma(1+ait)$ for any real $a$ and since $\Gamma(1) = 1$, we obtain  

\noindent{$(6.7)$}
\begin{align}
	C ({\rm sym}^{2} (\pi)) &= L(1, {\rm sym}^{2} (\pi_{\infty}) \times {\rm sym}^{2} (\pi_{\infty}))
	\notag \\
	&= \pi^{- 3} {|\Gamma(\frac{1}{2} + 2 i t)|}^{2} 
	{|\Gamma (\frac{1}{2} + i t)|}^{4}
	\notag \\
	&= 1 / (\cosh (2 \pi t) {\cosh (\pi t)}^{2})
	\notag
\end{align}

\bigskip

Recall that Theorem B proves that if $\pi$ is not of {\it solvable polyhedral type}, then $L(s, {\rm sym}^2(\pi) \times {\rm sym}^2(\pi))$ admits no Landau-Siegel zero. To put this to use we need the following

\medskip

\noindent{\bf Proposition 6.8} \, \it
	If $\pi$ is a spherical cusipdal representation on $GL(2)/\Q$, then
	$\pi$ is not of solvable polyhedral type.
\rm

\medskip

{\it Proof of Proposition} \,
At each prime $p$ (resp. $\infty$) let $\sigma_p$ (resp. $\sigma_\infty$) denote the $2$-dimensional representation of $W'_{\Q_p}$ (resp. $W_\R$) associated to $\pi_p$ (resp. $\pi_\infty$) by the local Langlands correspondence. By the naturality of this correspondence, we know that the conductors of $\pi_p$ and $\sigma_p$ agree at every $p$. 

On the other hand, as  $\pi$ is spherical, the conductor of $\pi$, which is the product of the conductors of all the $\pi_p$, is trivial. This implies that for every $p$, the conductor of $\sigma_p$, and hence also that of sym$^j(\sigma_p)$ is trivial for any $j \geq 1$. Appealing to the local correspondence again, we see that

\noindent{$(6.9)$} \it

For any $j \leq 4$, the automorphic representation ${\rm sym}^{j} (\pi)$ is spherical.
\rm

\medskip

For the definition of conductors, for any local field $k$, of representations of GL$(n, k)$ admitting a Whittaker model,  see [JPSS79].

\medskip

Assume that $\pi$ is of solvable polyhedral type, i.e., 
it is either dihedral or tetrahedral or octahedral.

Recall that if $\pi$ is {\it dihedral}, then $\pi$ is automorphically induced. i.e.
there exists an idele class character $\chi$ of a quadratic field $K$
s.t.
$$
\pi \, \simeq \, I^{\Q}_{K} (\chi).
\leqno(6.10)
$$

If $\pi$ is {\it tetrahedral}, then by [KSh2000], sym$^2(\pi)$ is cuspidal and moreover,
$$
{\rm sym}^{2} (\pi) \, \simeq \, I_{K}^{\Q} (\chi),
\leqno(6.11)
$$
for some idele class character $\chi$ of some cyclic extension
$K$ of degree $3$ over $\Q$. 

If $\pi$ is {\it octahedral}, then by [KSh2001], sym$^3(\pi)$ is cuspidal and
$$
{\rm sym}^3(\pi) \, \simeq \, {\rm sym}^3(\pi) \otimes \mu,
$$
for some quadratic Dirichlet character $\mu$. This implies, by [AC89], that
$$
{\rm sym}^{3} (\pi) \, \simeq \, I_{K}^{\Q} (\eta),
\leqno(6.12)
$$
for some cuspidal automorphic 
representaton $\eta$ of GL$(2, \A_K)$, with $K$ being the quadratic field
associated to $\mu$. 

\medskip

In view of (6.9), it suffices to show that some sym$^j(\pi)$ must be ramified, thus giving a contradiction.
Thanks to (6.10), (6.11) and (6.12), one is reduced to proving the following

\medskip

\noindent{\bf Lemma 6.13} \, \it Let $K/\Q$ be a cyclic extension of degree $\ell$, a prime, and let $\eta$ be a cuspidal automorphic representation of GL$(m, \A_K)$, $m \geq 1$. Then $I_K^\Q(\eta)$  is ramified at some $p$.
\rm

\medskip

One can be much more precise than this, but this crude statement is sufficient for our purposes.
However it should be noted that there are polyhedral rpresentations, for example of holomorphic type of weight $1$ for $F = \Q$, with prime conductor.

\medskip

{\it Proof of Lemma 6.13} \,  Put $\Pi = I_K^\Q(\eta)$. Since $\Q$ has class number $1$, $K/\Q$ is ramified. So there exists some prime $p$, and a place $u$ of $K$ above $p$,  such that $K_u/\Q_p$ is ramified of degree $\ell$. The local component $\Pi_p$ is simply $I_{K_u}^{\Q_p}(\eta_u)$, and it is enough to check that $\Pi_p$ must be ramified. If $\sigma_u$ is the $m$-dimensional representation of $W'_{K_u}$, then the conductor of $\Pi_p$ is the same as that of Ind$_{K_u}^{\Q_p}(\sigma_u)$. Moreover, $\sigma_u$ is semisimple and its conductor is divisible by that of Ind$_{K_u}^{\Q_p}(\sigma'_u)$ for any irreducible subrespresentation $\sigma'_u$ of $W'_{K_u}$. So it suffices to prove the following

\medskip

\noindent{\bf Sublemma 6.14} \, \it Let $E/F$ be a cyclic ramified extension of non-archimedean local fields, and let $\sigma$ be an irreducible $m$-dimensional representation of $W'_E$. Then Ind$_E^F(\tau)$ is ramified.
\rm

\medskip

{\it Proof of Sublemma} \,  Since $W'_{E}$ is $W_E \times {\rm SL}(2, \C)$, the irreducibility hypothesis implies that 
$$
\sigma \, \simeq \, \tau \otimes {\rm sym}^j(st),
$$
for some irreducible $\tau$ of $W_E$ and $j \geq 0$, where $st$ denotes the natural $2$-dimensional represenation of SL$(2, \C)$. Then
$$
{\rm Ind}_E^F(\sigma) \, \simeq \, {\rm Ind}_{W_E}^{W_F}(\tau) \otimes {\rm sym}^j(st).
\leqno(6.14)
$$

It suffices to prove that ${\rm Ind}_{W_E}^{W_F}(\tau)$ is ramified. Recall that there is a short exact sequence
$$
1 \, \rightarrow \, I_F \, \rightarrow W_F \, \rightarrow \, \Z \, \rightarrow \, 1,
\leqno(6.15)
$$
where $I_F$ denotes the inertia subgroup of Gal$(\overline F/F)$. If $\F_q$ is the residue field of $F$ and $\varphi$ the Frobenius $x \to x^q$, then $W_F$ is just the inverse image of the group of integral powers of the $\varphi$ under the natural map
$$
{\rm Gal}(\overline F/F) \, \rightarrow \, {\rm Gal}(\overline \F_q/\F_q).
$$

Suppose $\rho: = $Ind$_{W_E}^{W_F}(\tau)$ is unramified. Then by definition  $I_F$ must act trivially, and since the quotient $W_F/I_F$ is abelian,  $\rho$ is forced to be a sum of one dimensional, unramified representations. For this one must have
\begin{itemize}
\item[(i)] dim$(\tau) \, = \, 1$; \, and
\item[(ii)] $\tau^\theta \, \simeq \, \tau$, with $\theta$ denoting the non-trivial automorphism of $E/F$.
\end{itemize}
Consequently we have
$$
\rho \, \simeq \, \oplus_{i=0}^{[E:F]-1} \nu\delta^i,
\leqno(6.16)
$$
where $\nu$ is a character of $W_F$ extending $\tau$ and $\delta$ the character of $W_F$ associaated to $E/F$.But whatever $nu$ is, $\nu\delta^i$ will necessarily be ramified for some $i$ between $0$ and $[E:F]-1$. Thus $\rho =
$Ind$_{W_E}^{W_F}$ is ramified, contradicting the supposition that it is unramified. Done.

\qed

We have now proved Proposition 6.8.

\bigskip

Next we need the following two lemmas.

\medskip

\noindent{\bf Lemma 6.17} \, \it 
	Let $L(s) = \Sigma^{\infty}_{n = 1} \frac{b(n)}{n^{s}}$
	be an $L$-series with nonnegative coefficients, with $b(1) = 1$.
	Assume that $L(s)$ converges for $\Re s > 1$ with an analytic
	continuation to $\Re s > 0$.

	Let $M > 1$.
	Suppose that $L(s)$ satisfies the growth condition below on
	the line $\Re s = \frac{1}{2}$
	\[
		|L(\frac{1}{2} + i \gamma)| \le M {(|\gamma| + 1)}^{B}
	\]
	for some postive constant $B$.

	If $L(s)$ has no real zeros in the range
	\[
		1 - \frac{1}{log M} < s < 1
	\]
	then there exists an effective constant $c = c(B)$
	such that
	\[
		\text{\rm Res}_{s = 1} L(s) \ge \frac{c}{\log M}
	\]
\rm

\medskip

For a proof, see [GHLL94].

\medskip

\noindent{\bf Lemma 6.18} \, \it
	Let $L(s) = L(s, {\rm sym}^{2}(\pi) \times {\rm sym}^{2}(\pi))$.
	Then there exist absolute constants $A$ and $B$
	such that
	\[
		L(\frac{1}{2} + i \gamma) \le
		{(\lambda + 1)}^{A} {(|\gamma| + 1)}^{B}
	\]
\rm

\medskip

{\it Proof}. 

Note that, for any prime $p$, as $\pi_p$ is unramified, the $p$-part of $L(s)$
is the reciprocal of a polynomial in $p^{-s}$ of degree
$9$. Let $\alpha_{p}$, $\beta_{p}$ be the coefficients of the
Satake representation of $\pi_{p}$. Note that we assume that
$\pi$ is self-dual,  thus
$$
	{L_{p} (s)}^{-1} =
	(1 - \alpha_{p}^{4} p^{-s}) (1 - \beta_{p}^{4} p^{-s})
	{(1 - \alpha_{p}^{2} p^{-s})}^{2} {(1 -\beta_{p}^{2} p^{-s})}^{2}
	{(1 - p^{-s})}^{3}
\leqno(6.19)
$$
Now apply classical bound $\vert\alpha_p\vert < p^{1/4}$,  $\vert\alpha_p\vert < p^{1/4}$
on the coefficients; we know a much stronger bound now (cf. [K], Appendix 2), but the $1/4$ bound suffices for us. 
Then $L(s)$ is bounded by an absolute constant on
the line $\Re (s) = 2$. Also, $L(s)$ satisfies a functional
equation relating $s$ and $1 - s$. Thus, we get a bound for $L(s)$
on the line $\Re (s) = -1$. We claim that
the ratio of gamma factors arising from the
functional equation is bounded by a certain fixed power
of $\lambda$ and the imaginary part of $s$. 
In fact, the constants giving
the infinite type  of $L(s)$ are all imaginary as $\lambda > 50$ for Maass forms of kevel $1$. (The constants are real or purely imaginary, and the latter happens iff $\lambda \geq 1/4$, which is a difficult open problem for Maass forms of higher level.). Moreover, the self-duality of $\pi$ implies that the constant set is
symmetric about the real axis. So,
the norm of the ratio of the gamma factors is a product
of a constant and some terms of
$|\frac{\Gamma(1 + i t)}{\Gamma(- \frac{1}{2} - i t)}|$
where $t$ involves the imaginary part of $s$ and the 
constants (of the infinity type). Note that
$$
	\left| \frac{\Gamma (1 + i t)}{\Gamma (- \frac{1}{2} - i t)} \right|
	= \left| \frac{\Gamma (1 + i t)}{\Gamma (-\frac{1}{2} + i t)} \right|
	= {|t|}^{\frac{3}{2}} (1 + O (t^{-1})),
\leqno(6.20)
$$
since for $a \le \sigma \le b$ we have the estimation.
$$
	|\Gamma (s)| = \sqrt{2 \pi} e^{- \frac{\pi}{2}} {|t|}^{\sigma - 1/2}
	(1 + O (t^{-1}))
\leqno(6.21)
$$ 
where the implied constant depends only on $a$ and $b$. 
Hence the claim.

\medskip

As $\pi$
is spherical, so is ${\rm sym}^{4} (\pi)$. Hence we get
$$
	L(- 1 + i \gamma) << {(\lambda + 1)}^{A} {(|\gamma| + 1)}^{B}
\leqno(6.22)
$$
for certain constants $A$ and $B$. Applying the Phragmen-Lindel\"{o}f
principle in the strip $- 1 \le \Re (s) \le 2$, we see that
the same bound applies also on the line $\Re (s) = \frac{1}{2}$.

\qedsymbol

\medskip

The following proposition sets up the relationship between
the Peterson norm of the normalized automorphic 
function ${\rm sym}^{2}(g)$ and the residue of a certain
$L$-series at $s = 1$. 

\medskip

Denote $Z_{n} (\A)$ the center of $GL(n, \A)$. 
Denote $E^{*} (g, h_{s})$ the Eisenstein series, where
$h_{s} = \prod_{v} h_{s, v}$ and $h_{s, v}$
is in the space of the induced representation
$$
	{\rm Ind}_{P(n-1, 1, F_{v})}^{GL(n, F_{v})} (\delta_P^{s})
\leqno(6.23)
$$
where $\delta_P$ is the modular quasicharacter of
the standard parabolic subgroup $P(\Q_{v})$ of type $(n-1,1)$, whose
whose Levi facor is GL$(n - 1) \times {\rm GL}(1)$.

\medskip 

\noindent{\bf Proposition 6.24} \, \it
	Let $\Pi = \Pi_{\infty} \otimes \Pi_{f}$
	be an unramified cusp form on $GL(n, \A)$,
	with $\Pi_{\infty}$ a spherical principal
	series representation with trivial central character.
	Then   
$$
		\int_{Z_n(\A) GL(n, \Q) \backslash GL(n, \A)}
		\phi (g) {\overline \phi}) (g) E^{*} (g, h_{s}) dg =
		\frac{L (s, \Pi_{\infty} \times \Pi_{\infty}) L(s, \Pi_f \times \Pi_f)}
		{L (1, \Pi_{\infty} \times \Pi_{\infty})},
$$
where $\phi$ is the normalized function in 
the space of $\Pi$.
Furthermore, 
$$
		\langle\,\phi, \phi\,\rangle {\rm Res}_{s = 1} E^{*} (g, h_{s})
		\, = \, {\rm Res}_{s = 1} L(s, \Pi \times \Pi)
$$
\rm

\medskip

{\it Proof of Propsition 6.24}

Let us study the integral
$$
	I = \int_{Z_{n}(\A) GL(n, \Q) \backslash GL(n, \A)}
	\phi (g) \phi (g) E^{*} (g, h_{s}) dg.
\leqno(6.25)
$$
By the Rankin-Selberg unfolding method, we have
$$
	I = \prod_{v} \Psi(v, W_{\phi, v}, W_{\phi, v}, h_{s, v})
\leqno(6.26)
$$
where
$$
	I_{v} = \Psi(v, W, W, h_{s, v}) =
	\int_{Z_{n} (\Q_{v}) X_{n} (\Q_{v}) \backslash GL_{n} (\Q_{v})}
	W (g) W (g) h_{s} (g) dg
\leqno(6.27)
$$
Here $X_{n}$ denotes the subgroup of the upper triangular,
unipotent matrices, and $W_{\phi, v}$ is a Whittaker function
for $\Pi_{v}$. By Jacquet-Shalika [JS81], when we choose $\phi$ to be
a new vector,  this local integral $I_{v}$
equals to $L(s, \Pi_{v} \times \Pi_{v})$ when $v$ is nonarchimedean.

\medskip

When $v$ is archimedean, we appeal to the work of Stade ([St93], [St2001]) and obtain 
$$
I_{v} \, = \, {C (\Pi)}^{-1} L(s, \Pi_{v} \times \Pi_{v})
,\leqno(6.28)
$$
where $C (\Pi) = L(1, \Pi_{\infty} \times \Pi_{\infty})$. It appears that such a result has also been obtained by Jacquet and Shalika in the spherical case. In the non-spherical case, they can prove only that the $L$-factor is a fiite linear combination of such integrals.

\medskip

Thus $I$ is in fact the same as the quotient of the complete $L$-series for $\Pi \times \Pi$
by $L(1, \Pi_{\infty} \times \Pi_{\infty})$. Note that the Whitaker function
at infinity we take here differs from the standard one used by Stade in [St2001]
by the factor $C (\Pi)^{-1/2}$.

Now take the residue at $s=1$ on both sides, and note that
${\rm Res}_{s = 1} E^{*}(g, h_{s})$ is a positive
constant independent of $g$. Hence the Proposition.

\qedsymbol

\bigskip

{\it Proof of Corollary C} (contd).

Since sym$^2(\pi)$ is spherical in our case, we may apply
Proposition 6.24 with $\Pi = {\rm sym}^{2}(\pi)$
and get
$$
	(\,{\rm sym}^{2} (f), {\rm sym}^{2} (f)\,)
	= {C}^{-1} {\rm Res}_{s = 1} L(s, {\rm sym}^{2} (\pi) \times {\rm sym}^{2} (\pi))
\leqno(6.29)
$$where $C = {\rm Res}_{s = 1} E^{*}(g, h_{s})$.

\medskip

The right side of the corollary is easy since
$$
	{\rm Res}_{s = 1} L(s, {\rm sym}^{2} (\pi) \times {\rm sym}^{2} (\pi))
	= L(1, {\rm sym}^{2} (\pi)) L(1, {\rm sym}^{4} (\pi))
\leqno(6.30)
$$
which is bounded by any arbituary power of $1 + \lambda$.
(See [HRa95])

\medskip

To prove the left side, it suffices to show that
$$
	{\rm Res}_{s = 1} L(s, {\rm sym}^{2} (\pi) \times {\rm sym}^{2} (\pi))
	>> {(\log (1 + \lambda))}^{-1}
\leqno(6.31)
$$

For this apply Lemmas 6.17 and 6.18, with
$L(s) = L(s, {\rm sym}^{2} (\pi) \times {\rm sym}^{2} (\pi))$ and $M =
{(\lambda + 1)}^{B'}$ for suitably large constant $B'$. 

\bigskip

It remains to prove the asserted bound on the first Fourier coefficient of 
the spectral normalization of sym$^2(g)$. 
Put
$$
\Gamma \, = \, {\rm GL}(3, \Z),
\leqno(6.32)
$$
$$
\Gamma_0 \, = \, \{ \gamma = (\gamma_{ij}) \in \Gamma \, \vert \, \gamma_{31} = \gamma_{32} = 0, \, \gamma_{33} = 1\},
$$
and
$$
\Gamma_\infty \, = \, \{ \gamma = (\gamma_{ij}) \in \Gamma \, \vert \, \gamma_{ij} = 0 \, {\rm if} \, i > j\}.
$$

Recall that the cusp form sym$^2(g)$ on GL$(3)/\Q$ being of spherical type defines, and is determined by, a function, again denoted by sym$^2(g)$, on
the double coset space
$$
\Gamma\backslash {\rm GL}(3, \R)/Z_\R{\rm O}(3),
\leqno(6.33)
$$
where $Z_\R$ is the center of GL$(3, \R)$.

\medskip

Define the {\it spectrallly normalized} function in the space of sym$^2(\pi)$ to be
$$
{\rm sym}^2(g)^{\rm spec} \, = \, {\rm sym}^2(g)/||{\rm sym}^2(g)||,
\leqno(6.34)
$$
with $|| \, . \, ||$ denoting (as usual) the $L^2$-norm given by $\langle \, . \, , \, . \, \rangle^{1/2}$.

\medskip

The adelic Fourier expansion (6.4) gives rise to the following explicit expansion 
(see [Bu89], page 71, formula (2.1.6)) as a function of GL$(3, \R)$:

$$
{\rm sym}^2(g)(x) \, = \, \sum\limits_{(m,n) \ne (0,0)} \sum\limits_{\gamma \in \Gamma_\infty\backslash \Gamma_0} \, \frac{a(m,n)({\rm sym}^2(g))}{mn} W_\infty\left(
\begin{pmatrix}
mn & 0 & 0\\
0 & n & 0 \\
0 & 0 & 1
\end{pmatrix}
\gamma x\right).
\leqno(6.35)
$$
The coefficients $a(m,n)({\rm sym}^2(g))$ are bimultiplicative, implying in particular that the first coefficient $a(1,1)({\rm sym}^2(g))$ is equal to $1$. Consequently,
$$
a(1,1) : \, = \, a(1,1)({\rm sym}^2(g)^{\rm spec}) \, = \, \frac{1}{||{\rm sym}^2(g)||},
\leqno(6.36)
$$  
and
$$
|a(1,1)|^2\langle {\rm sym}^2(g)\, , \, {\rm sym}^2(g)\rangle \, = \, 1.
$$
Hence the bound on $|a(1,1)|$ follows from the bound proved above for
$\langle {\rm sym}^2(g)\, , \, {\rm sym}^2(g)\rangle$. Done.

\qedsymbol

\vskip 0.2in

\section{\bf Proof of  Theorem D}

\bigskip

Let $\pi$ be a cuspidal automorphic representation of GL$(2, \A_F)$ 
of trivial central character. Denote by $S$ the union of the set $S_\infty$ of archimedean places of $F$ with the set of finite places where $\pi$ is ramified. Given any Euler product $L(s) = \prod_v L_v(s)$ over $F$, we will write $L^S(s)$ to mean the (incomplete Euler) product of $L_v(s)$ over all $v$ outside $S$. 

Next recall (see section 3) that for every $j \leq 4$, there is an isobaric auutomorphic representation sym$^j(\pi)$ of GL$(j+1, \A_F)$, established long ago for $j=2$ by S. Gelbart and H. Jacquet [GJ77], and very recently for $j =3$, resp. $j=4$, by H. Kim and F. Shahidi ([KSh2000]), resp. H. Kim ([K2000]), such that
$$
L(s, {\rm sym}^j(\pi)) \, = \, L(s, \pi, {\rm sym}^j).
$$

\medskip

\noindent{\bf Lemma 7.1} \, \it Let $T$ be any finite set of places. Then we have the following factorizations of incomplete $L$-functions:
$$
L^T(s, {\rm sym}^3(\pi), {\rm sym}^2) \, = \, L^T(s, \pi, {\rm sym}^6)L^T(s, {\rm sym}^2(\pi))
\leqno(i)
$$
and
$$
L^T(s, {\rm sym}^4(\pi), {\rm sym}^2) \, = \, L^T(s, \pi, {\rm sym}^8)L^T(s, {\rm sym}^4(\pi))\zeta_F^T(s).
\leqno(ii)
$$
\rm

\medskip

{\it Proof} \, It suffices to prove these locally at every place otside $T$. But at any $v$, we have by definition,
$$
L(s, {\rm sym}^3(\pi_v), {\rm sym}^2) \, = \,  L(s, \Lambda^2({\rm sym}^4(\sigma_v))),
\leqno(7.2)
$$
and
$$
L(s, {\rm sym}^4(\pi_v), {\rm sym}^2) \, = \, L(s, {\rm sym}^2({\rm sym}^4(\sigma_v))),
$$
where $\sigma_v$ is the $2$-dimensional representation of $W_{F_v}$, resp. $W'_{F_v}$, associated to $\pi_v$ by the local correspondence for $v$ archimedean, resp. non-archimedean. By the Clebsch-Gordon identities, we have
$$
{\rm sym}^2({\rm sym}^3(\sigma_v)) \, \simeq \, {\rm sym}^6(\sigma_v) \oplus {\rm sym}^2(\sigma_v),
\leqno(7.3)
$$
and
$$
 {\rm sym}^2({\rm sym}^4(\sigma_v)) \, \simeq \, {\rm sym}^8(\sigma_v) \oplus {\rm sym}^4(\sigma_v) \oplus 1.
$$
The assertion of the Lemma now follows.

\qed

\medskip

\noindent{\bf Lemma 7.4} \, \it Let $\pi$ be a cuspidal automorphic representation of GL$(2, \A_F)$ with trivial central character, and let $v$ be a place where $\pi_v$ is a ramified, non-tempered principal series representatin. Then ${\rm sym}^4(\pi)$ is unramified at $v$.
\rm

\medskip

{\it Proof}. \, As $\pi_v$ is a ramified principal series representation of trivial central character, we must have
$$
\pi_v \, \simeq \, \mu_v \boxplus \mu_v^{-1},
\leqno(7.5)
$$
for a ramified (quasi-)character $\mu_v$ of $F_v^\ast$.  Since $\pi_v$ is non-tempered,
we may write, after possibly interchnging $\mu_v$ and $\mu_v^{-1}$,
$$
\mu_v \, = \, \nu_v\vert . \vert_v^t,
\leqno(7.6)
$$
for a unitary character $\nu_v$ of $F_v^\ast$ and a real number $t > 0$. ($\vert . \vert_v$ denotes as usual the normalizeed absolute value on $F_v$.) On the other hand, the unitarity of $\pi_v$ says that its complex conjugate representation $\overline \pi_v$ is isomorphic to the contragredient $\pi_v^\vee$. This forces the identity
$$
\nu_v \, = \, \overline \nu_v.
$$
Since $\nu_v$ is unitary, we get
$$
\nu_v^2 \, = \, 1 \quad {\rm and} \quad \pi_v \, \simeq \, \nu_v \otimes \pi^0_v,
\leqno(7.7)
$$
where
$$
\pi_v^0 \, \simeq \, \vert . \vert_v^t \boxplus \vert . \vert_v^{-1}.
$$
Then the associated $2$-dimensional Weil group representation $\sigma_v$ is of the form $\nu_v \otimes \sigma_v^0$, with $\sigma_v^0$ corresponding to $\pi_v^0$. Moreover, since $\nu_v$ is quadratic, we have for any $j \geq 1$,
$$
{\rm sym}^{2j}(\sigma_v) \, \simeq \,  {\rm sym}^{2j}(\sigma_v^0),
\leqno(7.8)
$$
which is unramified.
Since by [K2000], sym$^4(\pi)_v$ corresponds to sym$^4(\sigma_v)$ (at every place $v$), we see that it must be unramified as claimed.

\qed

\medskip

\noindent{\bf Lemma 7.9} \, \it Let $\pi$ be a cuspidal automorphic representation of GL$(2, \A_F)$, and $v$ a place of $F$ where $\pi_v$ is tempered.
Then for any $j \geq 1$, the local factor $L(s, \pi_v, {\rm sym}^j)$ is holomorphic in $\Re(s) > 1/2$ except for a possible pole at $s=1$..
\rm

\medskip

{\it Proof}. \, If $v$ is archimedean, or if $v$ is finite  but $\pi_v$ is not special, $\pi_v$ corresponds to a $2$-dimensional representation $\sigma_v$ of the local Weil group $W_{F_v}$. The temperedness of $\pi_v$ implies that $\sigma_v$ has bounded image in GL$(2, \C)$. Then for any finite-dimensional $\C$-representation $r$ of dimension $N$, in particular for sym$^j$, of GL$(2, \C)$, the image of $r(\sigma_v)$ will be bounded, and this implies that the admissible, irreducible representation $\Pi_v$ of GL$(N, F_v)$, associated to $r(\sigma_v)$ by the local Langlands correspondence, is tempered. Then $L(s,  \Pi_v)$ is holomorphic in $\Re(s) > 1/2$ except for a possible pole at $s=1$ (see [BaR94]). (One can also prove directly, using the extension in [De73] of Brauer's theorem to the representations of $W_{F_v}$, that $L(s, r(\sigma_v)$ has the requisite property.) We are now done in this case because $L(s, \pi_v; {\rm sym}^j)$ is defined to be $L(s, {\rm sym}^j(\sigma_v))$.

So we may take $v$ to be finite and assume that $\pi_v$ is a special representaton $sp(\mu_v)$ (see [HRa95], setion 2 for notation), associated to the partition $2 = 1 + 1$ and a (unitary) character $\mu_v$ of $F_v^\ast$. Then the associated $\sigma_v$ is of the form $(w,g) \to \mu_v(w) \otimes g$, for all $w$ in $W_{F_v}$ and $g$ in SL$(2, \C)$. So we have
$$
{\rm sym}^j(\sigma_v) \, \simeq \, \mu_v^j \otimes {\rm sym}^j,
$$
which corresponds to the special representation $sp(\mu_v^j)$ of GL$(j+1, F_v)$ associated to the partition $j+1 = 1 + \ldots + 1$ and the character $\mu_v^j$. Now we may appeal to the fact (see [BaR94]) that for any unitary character $\nu_v$, the function $L(s, sp(\nu_v))$ is holomorphic in $\Re(s) > 0$.

\qed

\medskip

Having established these preliminary lemmas, we are ready to begin the {\it proof of Theorem D}. Let $S$ denote the union of the archimedean places of $F$ with the set of finite places $v$ where $\pi_v$ is ramified {\it and} tempered. 
In view of Lemma 7.9, it suffices to show the following

\medskip

\noindent{\bf Proposition 7.10} \, \it The incomplete $L$-function $L^S(s, \pi; {\rm sym}^6)$ is holomorphic in the real interval $(1-\frac{c}{\log M}, 1)$ for a positive, effective constant $c$ independent of $\pi$, with $M$ denoting the thickened conductor of $\pi$. The same result holds for the symmetric $8$th power $L$-function if $F$ is a Galois extension of $\Q$ not containing any quadratic extension of $\Q$.
\rm

{\it Proof}. \, When $\pi$ is of solvable polyhedral type, the results of Kim and Shahidi in [Ksh2001]  imply that $L^S(s, \pi; {\rm sym}^j)$ is holomorphic in $(1/2,1)$ for any $j \leq 9$. So we may assume that we are not in this case, so that sym$^4(\pi)$ is a cuspidal automorphic representation of GL$(5, \A_F)$. 

By the definition of $S$, given any place $v$ outside $S$, $\pi_v$ is either unramified or a ramified, non-tempered principal series representation. Thanks to Lemma 7.4, sym$^4(\pi_v)$ is unramified in either case.  So we may appeal to the work of Bump-Ginzburg ([BuG92]) on the symmetric square $L$-fnctions, we get the holomorphy in $(1/2, 1)$ of the incomplete $L$-functions $L^S(s, {\rm sym}^4(\pi); \Lambda^2)$ and $L^S(s, {\rm sym}^4(\pi); {\rm sym}^2)$.

Next we appeal to the identities of Lemma 7.1 with $T = S$. The assertion of the Proposition is then clear for the symmetric $6$-th power $L$-function since $L^S(s, {\rm sym}^2(\pi))$ admits no Landau-Siegel zero by [GHLL94]. So let us turn our attention to the (incomplete) symmetric $8$-th power $L$-function of $\pi$. It suffices, by the identity (ii) of Lemma 7.1, that it suffices to show that $L^S(s, {\rm sym}^4(\pi))\zeta^S_F(s)$ admits no Landau-Siegel zero. Since $F$ is by hypothesis a Galois extension of $\Q$ not containing any quadratic field, one knows by Stark ([St]) that $\zeta_F^S(s)$ admits no Landau-Siegel zero. So we are finally done by our proof of Theorem B, where we showed that $L^S(s, {\rm sym}^4(\pi))$ admits no Landau-Siegel zero. Strictly speaking, we showed it for the full $L$-function. But the local factors at $S$, being tempered, do not have any pole in $(1/2,1)$. 

\qed

\vskip 0.2in

\section*{\bf Bibliography}

\begin{description}

\item[{[AC89]}] J. Arthur and L. Clozel, {\it Simple Algebras, Base Change
and the Advanced Theory of the Trace Formula}, Ann. Math. Studies {\bf 120}
(1989), Princeton, NJ.

\item[{[Ba97]}] W.~Banks, \emph{Twisted symmetric-square $L$-functions and the nonexistence of Siegel zeros on ${\rm GL}(3)$}, Duke Math. J.
{\bf 87} (1997), no. 2, 343--353.

\item[{[BaR94]}] L.~Barthel and D.~Ramakrishnan, \emph{A non-vanishing result for twists of $L$-functions of GL$(n)$}, Duke Math. Journal {\bf 74}, no.3 (1994), 681-700.

\item[{[Bu89]}] D.~Bump, \emph{The Rankin-Selberg method: A survey}, in {\it Number theory, trace formulas and discrete groups} (Oslo, 1987), 49--109,
Academic Press, Boston, MA (1989).

\item[{[BuG92]}] D.~Bump and D.~Ginzburg, \emph{Symmetric square $L$-functions on
${\rm GL}(r)$}, Ann. of Math. (2) {\bf 136}, no. 1, 137--205 (1992).

\item[{[CoPS94]}]
J.~Cogdell and I.~Piatetski-Shapiro, \emph{Converse Theorems
for $GL_n$}, Publications Math. IHES {\bf 79} (1994),
157--214.

\item[{[De73]}] P.~Deligne, \emph{Les constantes des \'equations fonctionnelles des
fonctions $L$}, in {\it Modular functions of one variable} II, Springer
Lecture Notes {\bf 349} (1973), 501-597.

\item[{[Ge75]}] S. Gelbart, {\it Automorphic forms on adele groups}, Annals of Math. Studies {\bf 83} (1975), Princeton.

\item[{[GJ79]}] S.~Gelbart and H.~Jacquet, 
\emph{A relation between automorphic
representations of GL$(2)$ and GL$(3)$},
 Ann. Scient. \'Ec. Norm. Sup. (4)
{\bf 11} (1979), 471--542.

\item[{[GHLL94]}] D.~Goldfeld, J.~Hoffstein and D.~Liemann,
		\emph{An effective zero free region}, Ann. of Math. {\bf 140} (1994),
		appendix to [HL94].
		
\item[{[GS2000]}] A.~Granville and H.~Stark, \emph{$abc$ implies no ``Siegel zeros'' 
for $L$-functions of characters with negative discriminant}, Inventiones Math. 
{\bf 139} (2000), no. 3, 509--523.

\item[{[HaT2000]}] M.~Harris and R.~Taylor, \emph{On the geometry and cohomology of some 
simple Shimura varieties}, preprint (2000), to appear in the Annals of Math. Studies, Princeton.

\item[{[He2000]}] G.~Henniart, \emph{Une preuve simple des conjectures de Langlands pour 
${\rm GL}(n)$ sur un corps $p$-adique}, Invent. Math. {\bf 139}, no. 2, 439--455
(2000).

\item[{[HL94]}] J.~Hoffstein and P.~Lockhart,
		\emph{Coefficients of Mass forms and the Siegel zero},
		Ann. Math. (2) \textbf{140} (1994), 161--181. 

\item[{[HRa95]}] J.~Hoffstein and D.~Ramakrishnan,
		\emph{Siegel Zeros and Cusp Forms},
		IMRN(1995),  No.\ \textbf{6},  279--308.

\item[{[IwS2000]}] H.~Iwaniec and P.~Sarnak, \emph{The non-vanishing of central values 
of automorphic $L$-functions and Landau-Siegel zeros}, Israel Journal of
Math. 120 (2000), part A, 155--177.

\item[{[JPSS79]}] H.~Jacquet, I.~Piatetski-Shapiro and J.A.~Shalika,
\emph{Automorphic forms on ${\rm GL}(3)$}. II. Ann. of Math. (2) {\bf 109}, 
no. 2, 213--258 (1979).  

\item[{[JPSS83]}] H.~Jacquet, I.~Piatetski-Shapiro and J.A.~Shalika,
\emph{Rankin-Selberg convolutions}, Amer. J of Math. {\bf 105} (1983), 367--464.

\item[{[JS81]}] H.~Jacquet and J.A.~Shalika, 
\emph{Euler products and
the classification of automorphic forms} I \& II, Amer. J of
Math. {\bf 103} (1981), 499--558 \& 777--815.

\item[{[JS90]}] H.~Jacquet and J.A.~Shalika, \emph{Rankin-Selberg convolutions: archimedean theory}, in {\it Piatetski-Shapiro Festschrift}, Israel Math. conf. Proc., Part II, 125-207, The Weizmann Science Press of Israel (1990).

\item[{[K2000]}] H.~Kim, \emph{Functoriality of the exterior square of GL$_4$ and the
symmetric fourth of GL$_2$}, preprint (2000).

\item[{[KSh2000]}] H.~Kim and F.~Shahidi,
\emph{Functorial products for GL$(2) \times $GL$(3)$ and the symmetric cube for GL$(2)$}, preprint (2000), to appear in Annals of Math.

\item[{[KSh2001]}] H.~Kim and F.~Shahidi, \emph{Cuspidality of symmetric powers with applications}, preprint (2001), to appear in the Duke Journal of Math.

\item[{[La70]}] R.P.~Langlands, \emph{Problems in the theory of automorphic forms}, in
{\it Lectures in modern analysis and applications III}, 
Lecture Notes in Math. {\bf 170} (1970), Springer-Verlag, Berlin, 18--61.

\item[{[La79]}] R.P.~Langlands, \emph{On the notion of an automorphic representation. A supplement}, in {\it Automorphic forms, Representations and $L$-functions}, ed. by A. Borel and W. Casselman, Proc. symp. Pure Math {\bf 33}, part 1, 203-207, AMS. Providence (1979). 

\item[{[La80]}] R.P.~Langlands, {\it Base change for GL$(2)$}, Annals of Math. Studies {\bf 96}, Princeton (1980).

\item[{[MW89]}] C.~Moeglin and J.-L.~Waldspurger, \emph{Poles des fonctions $L$ de paires pour GL$(N)$}, Appendice, Ann. Sci. \'Ecole Norm. Sup. (4) {\bf 22}, 667-674 (1989).

\item[{[Mo85]}] C.J.~Moreno,
		\emph{Analytic proof of the strong multiplicity one theorem},
		Amer. J.Math. \textbf{107}, \textsl{no} 1, 163--206.

\item[{[Mu94]}] Murty,~M.~Ram,
\emph{Selberg's conjectures and Artin $L$-functions}, 
Bull. Amer. Math. Soc. (N.S.) {\bf 31} (1994), no. 1, 1--14.

\item[{[PPS89]}] S.~.J.~Patterson and I.~Piatetski-Shapiro, \emph{The symmetric-square $L$-function 
attached to a cuspidal automorphic representation of ${\rm GL}\sb 3$}, 
Math. Ann. {\bf 283} (1989), no. 4, 551--572. 

\item[{[Ra99]}] D.~Ramakrishnan,
		\emph{Landau--Siegel Zeros and Cusp Forms},
		IAS Lecture (1999), preprint on www.math.caltech.edu/people/dinakar.html.

\item[{[Ra2000]}] D.~Ramakrishnan, \emph{Modularity of the Rankin-Selberg $L$-series, and
Multiplicity one for SL$(2)$}, Annals of Mathematics {\bf 152} (2000), 45--111.  

\item[{[Sh88]}] F.~Shahidi, \emph{On the Ramanujan conjecture and the
finiteness of poles for certain $L$-functions}, Ann. of Math. (2) {\bf 127}
(1988), 547--584.

\item[{[Sh90]}] F.~Shahidi, \emph{A proof of the Langlands conjecture on
Plancherel measures; Complementary series for $p$-adic groups}, Ann. of
Math. {\bf 132} (1990), 273-330.

\item[{[St93]}] E.~Stade, \emph{Hypergeometric series and Euler factors at infinity for $L$-functions on GL$(3, \R) \times $GL$(3, \R)$}, American Journal of Math. {\bf 115}, No.2 (1993), 371--387. 

\item[{[St2001]}] E.~Stade,
\emph{Archimedean $L$-Factors on $GL(n) \times GL(n)$
	 and Generalized Barnes Integrals}, to appear.

\item[{[Stk74]}] H. Stark, Some effective cases of the Brauer-Siegel theorem, Inventiones Math. 
{\bf 23} (1974), 135--152.

\bigskip

\end{description}

\vskip 0.3in

Dinakar Ramakrishnan \qquad \qquad Song Wang

\bigskip

\vskip 0.2in

\end{document}